\documentclass[reqno,10pt]{amsart}

\usepackage{xcolor}
\usepackage{graphicx}
\usepackage[hidelinks]{hyperref}

\usepackage{amscd}

\usepackage{amsmath}
\usepackage{amsthm}
\usepackage{amssymb}
\usepackage{latexsym,array}
\usepackage{amsfonts}
\usepackage{shadow}
\usepackage{amsbsy}
\usepackage{amssymb}
\usepackage{dsfont}
\usepackage{mathtools}
\usepackage{enumitem}
\usepackage{doi}

\usepackage{tikz-cd}
\usepackage{color}
\usepackage{graphicx}
\usepackage[hidelinks]{hyperref}
\usepackage{cleveref}
\usepackage{fullpage}
\usepackage{comment}

\usepackage{dsfont}

\usepackage{cleveref}

\numberwithin{equation}{section}

\newtheorem{theorem}{Theorem}[section]
\newtheorem{lemma}[theorem]{Lemma}

\newtheorem{proposition}[theorem]{Proposition}

\theoremstyle{definition}
\newtheorem{remark}[theorem]{{\bf Remark}}
\newtheorem{definition}[theorem]{Definition}

\usepackage{amscd}
\usepackage{tikz-cd}
\usepackage{tikz}

\newtheorem{example}[theorem]{Example}

\crefname{enumi}{}{}
\crefname{enumii}{}{}

\title[]{The Clifford Short-Time Fourier Transform with Non-Radial Window Functions and Its Uncertainty Principles}

\author[A. De Martino]{Antonino De Martino}
\address{(ADM)
	Politecnico di Milano\\Dipartimento di Matematica\\Via E. Bonardi, 9\\20133
	Milano, Italy
} \email{antonino.demartino@polimi.it}

\author[Yi-Qiao Xu]{Yi-Qiao Xu}
\address{(YQX)
 Beijing Institute of Technology,
School of Mathematics and Statistics,
No. 5 Zhongguancun South Street, Haidian District,
100081 Beijing, China
}\email{yiqiaoxu1016@163.com}

\date{}

\begin{document}

	\maketitle

\begin{abstract}
We develop a theory of the Clifford short-time Fourier transform (CSTFT) with non-radial window functions in even dimensions. The main difficulty is that the generalized Clifford translation no longer coincides with the ordinary translation. To overcome this, we establish an angular estimate for the Clifford–Fourier kernel, which yields a uniform $L^2$-bound for generalized translates of Schwartz functions. We derive a weighted Moyal identity involving a Clifford-valued weight and introduce an admissibility condition on the window function that allows us to obtain reconstruction and reproducing kernel formulas. We also provide an explicit example of a non-radial admissible window. Finally, we establish Pitt-type, logarithmic and Heisenberg-type uncertainty principles for the CSTFT with non-radial windows. 
\end{abstract}



%\tableofcontents

\section{Introduction}

The short-time Fourier transform (STFT) is a fundamental tool in time-frequency analysis, providing a description of the local frequency content of a signal through translations and modulations of a suitable window function. Its classical theory is well established; see, for example, \cite{G}. Among its many applications, uncertainty principles play an important role and have been investigated from various perspectives; see \cite{BJ, C, Don, GK, G1, lieb, W} and the references therein.

Several generalizations and applications of the classical Fourier transform have been studied in different mathematical settings; see, for example, \cite{Dunkl,DeBieCFT,DeBieFourier,DeJ,T}. Among these, the Clifford-Fourier transform arises naturally in Clifford analysis, a higher-dimensional extension of classical harmonic analysis; see \cite{red, ColomboSabadiniStruppa2011, green}.
Unlike the classical Fourier transform, the Clifford-Fourier transform has a kernel that is
generally unbounded, with polynomial growth depending on the dimension. 
More precisely, the Clifford-Fourier transform is an integral transform acting on Clifford-valued functions, defined by
$$\mathcal F_\pm f(y)=(2\pi)^{-d/2}\int_{\mathbb R^d}K_\pm(t,y)f(t)\,dt.$$
For even $d\geq2$, its kernel satisfies the polynomial bound
$$ |K_\pm(t,y)| \leq C_d(1+|t|)^{\lambda} (1+|y|)^{\lambda}, \qquad \lambda=\frac{d-2}{2}. $$
Such estimates arise in the general framework of Fourier-like integral transforms with
polynomially bounded kernels, which also includes the Fourier-Dunkl transform; see \cite{DeJ, GhobberJaming}. Moreover, in this paper we establish the following angular estimate for the Clifford-Fourier kernel: for every pair of nonzero, non-parallel vectors $t,y\in\mathbb R^d$,
\begin{equation}\label{I1}
	|K_-(t,y)| \leq \frac{C_d}{(\sin\theta)^{\lambda}},
\end{equation}
where $\theta\in(0,\pi)$ denotes the angle between $t$ and $y$. A Clifford short-time Fourier transform (CSTFT) was introduced in \cite{DeMartinoCSTFT} using
the Clifford-Fourier transform and its associated generalized translation operator. In that work, the window functions were assumed to be radial. This assumption considerably simplifies the analysis, since the generalized translation coincides with the classical translation for
radial functions. For non-radial windows, this property is no longer available, and several
fundamental features of the classical STFT must be reconsidered.

The main purpose of this paper is to develop the theory of the CSTFT with non-radial window functions and to establish some corresponding uncertainty principles. We begin by investigating the generalized translation operator given by

$$(\tau_x f)(t)=(2\pi)^{-d/2}\int_{\mathbb R^d}K_-(\xi,t)K_-(x,\xi)\mathcal F_-f(\xi)\,d\xi,
	\qquad x,t\in\mathbb R^d.
$$
We prove that it preserves the Schwartz space and by using \eqref{I1} we establish the uniform estimate
$$
\sup_{t\in\mathbb R^d} \int_{\mathbb R^d}|(\tau_xg)(t)|^2\,dx <\infty
$$
for every $g\in\mathcal S(\mathbb R^d;\mathbb R_{0,d})$ and every even dimension $d\geq2$.  We then investigate the structural properties of the CSTFT with general windows. A central
role is played by the weight
$$
W_{g_1,g_2}(t)= \int_{\mathbb R^d} \overline{(\tau_xg_1)(t)} (\tau_xg_2)(t)\,dx.
$$
Unlike the classical setting, the weight $W_{g_1,g_2}(t)$ is generally Clifford-valued and depends on the variable $t$. Consequently, the weighted Moyal identity involves a noncommutative weight rather than a constant scalar factor. This difference is particularly
relevant for reconstruction, where the invertibility of the weight cannot be taken
for granted. To address this issue, we introduce an admissibility condition based on the
quadratic form
$$Q_g(t,u)=\int_{\mathbb R^d}|(\tau_xg)(t)u|^2\,dx.$$
More precisely, we require that
$$Q_g(t,u)\geq c_g(1+|t|)^{-N}|u|^2$$
for some $c_g>0$ and $N\geq0$. This condition guarantees the invertibility
of $W_g(t)$ and a polynomial bound for its inverse, allowing us to establish reconstruction
and reproducing kernel formulas. We also prove that admissibility is stable under sufficiently small Schwartz perturbations of nonzero real-valued radial windows. In particular, we construct explicit non-radial Clifford-valued admissible windows in every even dimension.

Finally, we establish some  uncertainty principles for the CSTFT with non-radial windows. Starting from Pitt's inequality for the Clifford-Fourier transform, we derive a corresponding inequality for the CSTFT. We subsequently obtain a logarithmic uncertainty principle and a Heisenberg-type inequality.

%\newline
%\newline
%\emph{The plan of the paper}
%\newline
%\newline
%The paper is organized as follows. In Section~2 we recall the basic notions of Clifford algebras and Clifford analysis. Section~3 recalls the main properties of the Clifford-Fourier transform and its associated generalized translation operator that will be used
%throughout the paper, and establishes some new important properties. In Section~4 we
%study the CSTFT with general window functions. Section~5 is devoted to uncertainty principles. We prove Pitt-type, logarithmic, Heisenberg-type and Lieb-type uncertainty inequalities and conclude the section with local uncertainty principles for the CSTFT.

\section{Preliminaries}

In this section, we recall the main notation and basic notions used throughout the paper. We begin by recalling some basic facts on Clifford algebras. Let $e_1,\ldots,e_d$ be an orthonormal basis of the Euclidean space $\mathbb{R}^d$, satisfying the relations
$e_je_k+e_ke_j=-2\delta_{jk}$, $1\leq j,k\leq d,$ where $\delta_{jk}$ denotes the Kronecker delta. The real Clifford algebra $\mathbb{R}_{0,d}$ is the real associative algebra generated by these elements. A basis of $\mathbb{R}_{0,d}$, viewed as a real vector space, is given by the elements $e_S=e_{j_1}\cdots e_{j_k},$ where $S=\{j_1,\ldots,j_k\}\subseteq\{1,\ldots,d\}$ with $1\leq j_1<\cdots<j_k\leq d$, together with $e_{\emptyset}=1$. Hence every element $a\in\mathbb{R}_{0,d}$ can be written as
$$
a=\sum_S a_S e_S,\qquad a_S\in\mathbb{R}.
$$
The Clifford conjugate of $a$ is defined by
$$
\overline{a}=\sum_S (-1)^{\frac{|S|(|S|+1)}{2}} a_S e_S.
$$
When complex scalars occur, this conjugation is extended anti-linearly, so that
\begin{equation}
\label{reverse}
\overline{\lambda}=\lambda^c,\qquad \overline{e_i}=-e_i,\quad i=1,\ldots,d,\qquad
 \overline{ab}=\overline b\,\overline a,
\end{equation}
where $\lambda^c$ denotes the complex conjugate of $\lambda$. We endow $\mathbb R_{0,d}$ with the standard real inner product
\begin{equation}
\label{inner1}
(a,b):=\operatorname{Sc}(\overline{a}b),
\end{equation}
whose associated norm is
$$
|a|^2=(a,a)=\operatorname{Sc}(\overline{a}a)=\left(\sum_S |a_S|^2\right).
$$
The norm on $\mathbb{R}_{0,d}$ satisfies the following inequality
\begin{equation}
	\label{estmou}
	|ab| \leq 2^{\frac{d}{2}}|a||b|, \qquad a,b\in\mathbb R_{0,d}.
\end{equation}
An element $x=(x_1,...,x_d) \in \mathbb{R}^d$ can be identified with the vector $x=\sum_{j=1}^d x_j e_j$. The multiplication of two vectors $x$ and $y$ is given by
$$xy=-\langle x, y \rangle+ x \wedge y,$$
where the inner and the wedge product are defined, respectively by
\begin{equation}
\label{product}
\langle x,y\rangle=\sum_{j=1}^d x_jy_j=-\frac{1}{2}(xy+yx), \qquad x \wedge y= \sum_{i<k}e_i e_k (x_i y_k-x_k y_i)=\frac{1}{2}(xy-yx).
\end{equation}

We denote the space $L^p(\mathbb{R}^d;\mathbb{R}_{0,d})$ as the module of Clifford-valued functions $f: \mathbb{R}^d \to \mathbb{R}_{0,d}$ with finite norm
\begin{equation}
	\|f\|_{L^p(\mathbb{R}^d;\mathbb{R}_{0,d})}=
	\begin{cases}
		\left(\displaystyle\int_{\mathbb R^d}
		|f(x)|^p\,dx\right)^{1/p},
		& 1 \leq p < \infty,\\[1ex]
		\operatorname*{ess\,sup}_{x \in \mathbb{R}^d}
		|f(x)|,
		& p=\infty.
	\end{cases}
\end{equation}
where $dx=dx_1\cdots dx_d$ represents the usual Lebesgue in $\mathbb{R}^d$. In particular, the norm on $L^2(\mathbb{R}^d;\mathbb{R}_{0,d})$ is induced by the scalar inner product
\begin{equation}
\langle f,g\rangle := \int_{\mathbb{R}^d} \operatorname{Sc}\bigl(\overline{f(x)}\,g(x)\bigr)\,dx.
\end{equation}

We denote by $\mathcal{S}(\mathbb{R}^d;\mathbb{R}_{0,d})$ the Schwartz space of Clifford-valued functions, defined by
\begin{eqnarray*}
 \mathcal S(\mathbb R^d;\mathbb{R}_{0,d})&=& \{f: \mathbb{R}^d \to \mathbb{R}_{0,d} \, : \, \sup_{x \in \mathbb{R}^d} | x^\alpha \partial^\beta f(x)|<\infty; \quad x^\alpha=x^{\alpha_1}_1x^{\alpha_2}_2...x^{\alpha_d}_d, \, \alpha=(\alpha_1,...,\alpha_d), \, \beta \in \mathbb{N}^d_0\}.
\end{eqnarray*}

\section{The Clifford Fourier transform}

In \cite{H} Howe found an exponential operator equivalent definition of the classical complex Fourier transform using the harmonic oscillator $\Delta-|x|^2$, i.e.
$$ \mathcal{F}= \hbox{exp} \left(- \frac{\pi i}{2} \mathcal{H}\right)=\sum_{n=0}^{\infty} \frac{1}{n!} \left(-\frac{i \pi}{2}\right)^n \mathcal{H}^n$$
where $ \mathcal{H}$ is the scalar-valued operator given by $\mathcal{H}=\frac{1}{2} (-\Delta+|x|^2-d)$, and $\Delta$ is the Laplace operator. The Clifford-Fourier transform was first introduced in \cite{BSS} by decomposing the operator $\mathcal{H}$ into a sum of Clifford algebra-valued second-order operators involving the Gamma operator $\Gamma_x$. More precisely, the Clifford-Fourier transform was defined as
$$
\mathcal F_\pm=e^{\frac{i\pi d}{4}}e^{\frac{i\pi}{4}(\Delta-|x|^2\mp2\Gamma_x)}=e^{\frac{i\pi d}{4}}e^{\mp\frac{i\pi}{2}\Gamma_x}e^{\frac{i\pi}{4}(\Delta-|x|^2)}.
$$
where $\Gamma_x$ is the Gamma operator associated to a vector $x$, given by
\begin{equation}
\label{gamma}
\Gamma_x=-\sum_{i<k} e_i e_k\bigl(x_i\partial_{x_k}-x_k\partial_{x_i}\bigr).
\end{equation}

In \cite{DeBieCFT}, the authors provided an integral representation of the above expression.

\begin{definition}[Clifford-Fourier transform (CFT)]
\label{cft}
On the class of functions $\mathcal S(\mathbb R^d;\mathbb{R}_{0,d})$, we define the Clifford-Fourier transform as
$$ \mathcal F_\pm f(y) =(2\pi)^{-d/2}\int_{\mathbb R^d}K_\pm(t,y)f(t)\,dt,$$
and their inverse as
\begin{equation}
\label{inverse}
\mathcal{F}_{\pm}^{-1}f(y) =(2\pi)^{-d/2} \int_{\mathbb R^d}\widetilde{K_{\pm}}(t,y) f(t)\,dt
\end{equation}
where
$$
K_{\pm}(t,y):=\exp\!\left(\mp\frac{i\pi}{2}\Gamma_y\right)e^{-i\langle t,y\rangle}, \quad \hbox{and} \quad \widetilde{K_{\pm}}(t,y):=\exp\!\left(\pm\frac{i\pi}{2}\Gamma_y\right)e^{i\langle t,y\rangle}.
	$$
\end{definition}

\begin{remark}
In the remainder of the paper, we will only consider the integral kernel $K_{-}(t,y)$. Indeed, the results for the other kernel can be obtained from the relation $K_{-}(t,y)=K_{+}(t,-y)$.
\end{remark}

\begin{remark}
We note that the kernels are not symmetric in the sense that
$K_{\pm}(t,y) \neq K_{\pm}(y,t)$, see \cite{DeBieCFT}. Hence, we adopt the convention that we always integrate with respect to the first variable of the kernel.
\end{remark}

A challenging problem has been to find an explicit expression for the kernel $K_{\pm}(t,y)$. In \cite{BSS1}, a two-dimensional kernel was obtained using techniques from Clifford analysis. More precisely, this is given by

\begin{equation}
\label{dim2}
K_{\pm}(t,y)=\cos w \mp e_1e_2\sin w, \qquad w=t_1y_2-t_2y_1.
\end{equation}

A complete expression for all even dimensions was given in \cite[Thm. 4.3]{DeBieCFT}, where the authors expressed the kernel $K_{\pm}(t,y)$ as a finite sum involving Bessel functions $J_{\alpha}(\rho)$. More precisely, it is given by
\begin{equation}
	K_{-}(t,y)=(-1)^{d/2}\left(\frac{\pi}{2}\right)^{1/2}\left(A_\lambda^*(s,\rho)+B_\lambda^*(s,\rho)+(t\wedge y)C_\lambda^*(s,\rho)\right), \qquad t,y \in \mathbb{R}^d,
	\label{kernelreal}
\end{equation}
where
\begin{equation}
	\label{eq1}
	A_\lambda^*(s,\rho)=\sum_{\ell=0}^{\left\lfloor\frac d4-\frac34\right\rfloor}
	\frac{\Gamma(d/2)}{2^\ell\ell!\,\Gamma(d/2-2\ell-1)}s^{d/2-2-2\ell}\widetilde J_{(d-2\ell-3)/2}(\rho)
\end{equation}

\begin{equation}
	\label{eq2}
	B_\lambda^*(s,\rho)=-\sum_{\ell=0}^{\left\lfloor\frac d4-\frac12\right\rfloor}\frac{\Gamma(d/2)}{2^\ell\ell!\,\Gamma(d/2-2\ell)}
	s^{d/2-1-2\ell}\widetilde J_{(d-2\ell-3)/2}(\rho),
\end{equation}

\begin{equation}
	\label{eq3}
	C_\lambda^*(s,\rho)=-\sum_{\ell=0}^{\left\lfloor\frac d4-\frac12\right\rfloor}
	\frac{\Gamma(d/2)}{2^\ell\ell!\,\Gamma(d/2-2\ell)}s^{d/2-1-2\ell}\widetilde J_{(d-2\ell-1)/2}(\rho),
\end{equation}
with $\lambda= \frac{d-2}{2}$, $s=\langle t,y\rangle$, $\rho=|t\wedge y| =\sqrt{|t|^2|y|^2-s^2}$, and $\widetilde J_\nu(\rho)=\rho^{-\nu}J_\nu(\rho)$.
\\In contrast with the classical kernel of the Fourier transform, the kernel $K_{-}(t,y)$ satisfies a polynomial bound, see \cite[Theorem~5.3]{DeBieCFT}.

\begin{proposition}
\label{est}
Let $d \geq 2$ be even and set $\lambda:= \frac{d-2}{2}$. For $t$, $y \in \mathbb{R}^d$, there exists a constant $C_d$ independent of $t$ and $y$ such that
$$ |K_{-}(t,y)| \le C_d(1+|t|)^{\lambda}(1+|y|)^{\lambda}.$$
\end{proposition}

By using Proposition \ref{est} in \cite[Thm. 6.1]{DeBieCFT} the authors specify the domain of the Clifford-Fourier transform.

\begin{theorem}
\label{cont1}
Let $d \geq 2$ be even. The Clifford-Fourier transform is well defined on $\mathcal{B}(\mathbb R^d;\mathbb{R}_{0,d})$, where
$$\mathcal B(\mathbb R^d;\mathbb{R}_{0,d})=\left\{
f\in L^1(\mathbb R^d;\mathbb{R}_{0,d}): \, \,
\|f\|_{\mathcal{B}}:=\int_{\mathbb R^d}(1+|x|)^\lambda
|f(x)|\,dx<\infty
\right\}.
$$
Moreover, the Clifford-Fourier transforms $ \mathcal{F}_{\pm}$ are continuous operators on $ \mathcal{S}(\mathbb{R}^d; \mathbb{R}_{0,d}).$
\end{theorem}

In \cite{LiLengFeiPW} the authors prove the following counterpart of the Plancherel and Parseval identity for the Clifford Fourier transform.

\begin{theorem}
	\label{par}
	If $f$, $g \in \mathcal{S}(\mathbb{R}^d; \mathbb{R}_{0,d})$ then we have
	$$ \int_{\mathbb{R}^d} \overline{f(y)}g(y)dy= \int_{\mathbb{R}^{d}} \overline{\mathcal{F}_{-}(f)(x)} \mathcal{F}_{-}(g)(x)dx.$$
	Moreover if $f=g$ we have $\| f\|_{2}=\|\mathcal{F}_{-} (f)\|_{2}$.
\end{theorem}

\begin{remark}
Since $\mathcal S(\mathbb R^d;\mathbb R_{0,d})$ is dense in $L^2(\mathbb R^d;\mathbb R_{0,d})$, the Plancherel identity in Theorem~\ref{par} allows us to extend $\mathcal F_-$ uniquely to an isometry on $L^2(\mathbb R^d;\mathbb R_{0,d})$. Moreover, by density and the Cauchy-Schwarz inequality, the Parseval identity also extends to all $f,g\in L^2(\mathbb R^d;\mathbb R_{0,d})$. These identities will also be used in their $L^2$ formulations when needed.
\end{remark}

Now, we prove a new estimate for the kernel that will be crucial in the rest of the paper.

\begin{lemma}\label{angular-kernel}
Let $d\geq2$ be even and set $\lambda=(d-2)/2$. For every pair of nonzero, non-parallel vectors $t,y\in\mathbb R^d$, we have
\begin{equation}\label{angularkernel}
	|K_-(t,y)| \leq \frac{C_d}{(\sin\theta)^\lambda},
\end{equation}
where $\theta\in(0,\pi)$ denotes the angle between $t$ and $y$, and $C_d>0$ depends only on the dimension.
\end{lemma}

\begin{proof}
For $d=2$, the result follows immediately from \eqref{dim2}. We estimate the component of the kernel in \eqref{kernelreal}. By the power series expansion near the origin and the asymptotic behavior of Bessel functions
at infinity, see \cite[formulas 9.1.7 and 9.2.5]{AS}, for every fixed $\nu\geq1/2$ there exists $C_\nu>0$ such that
\begin{equation}\label{Besselangularest}
|\widetilde J_\nu(\rho)| \leq C_\nu(1+\rho)^{-\nu-1/2},\qquad \rho\geq0.
\end{equation}
Let $\theta\in(0,\pi)$ be the angle between $t$ and $y$. Since $s=\langle t, y \rangle =|t||y|\cos\theta$ and $\rho=|t \wedge y|=|t||y|\sin\theta$ we have
\begin{equation}\label{scot}
|s|=\rho|\cot\theta|.
\end{equation}
We first estimate $A_\lambda^*(s,\rho)$. For $0\leq\ell\leq\lfloor d/4-3/4\rfloor$, set
$p:=\frac d2-2-2\ell$. By \eqref{Besselangularest} and \eqref{scot}, we obtain
$$|s|^p \left|\widetilde J_{(d-2\ell-3)/2} (\rho)\right| \leq C_{d,\ell}\frac{|s|^p}{(1+\rho)^{d/2-\ell-1}}=C_{d,\ell}|\cot\theta|^p \frac{\rho^p}{(1+\rho)^{p+\ell+1}}\leq \mathbf{c}_d|\cot\theta|^p,
$$
where $\mathbf{c}_d:= \max_{0\leq\ell\leq\lfloor d/4-3/4\rfloor}C_{d,\ell}$. Since $0\leq p\leq d/2-2$ and the sum is finite, we obtain
\begin{equation}\label{A-angular}
|A_\lambda^*(s,\rho)| \leq C_d(1+|\cot\theta|)^{\lambda-1}.
\end{equation}
Next, we consider $B_\lambda^*(s,\rho)$. For $0\leq\ell\leq\lfloor d/4-1/2\rfloor$, set
$ q:=\frac d2-1-2\ell$. Using again \eqref{Besselangularest} and \eqref{scot}, we get
$$|s|^q \left|\widetilde J_{(d-2\ell-3)/2}(\rho)\right|\leq
C_{d,\ell}\frac{|s|^q}{(1+\rho)^{d/2-\ell-1}}=C_{d,\ell}|\cot\theta|^q\frac{\rho^q}{(1+\rho)^{q+\ell}}\leq \mathbf{c}_d|\cot\theta|^q.
$$
Since $0\leq q\leq d/2-1=\lambda$, and the sum is finite we get
\begin{equation}\label{B-angular}
|B_\lambda^*(s,\rho)|\leq C_d(1+|\cot\theta|)^\lambda.
\end{equation}	
Finally, we estimate the bivector term. Since $|t\wedge y|=\rho$, by \eqref{Besselangularest} and \eqref{scot}, we obtain
$$\rho|s|^q\left|\widetilde J_{(d-2\ell-1)/2}(\rho)\right| \leq
C_{d, \ell}\frac{\rho|s|^q}{(1+\rho)^{d/2-\ell}}=C_{d, \ell}|\cot\theta|^q \frac{\rho^{q+1}}{(1+\rho)^{q+\ell+1}}\leq \mathbf{c}_d|\cot\theta|^q.
$$
Summing over the finitely many values of $\ell$, we obtain
\begin{equation}\label{C-angular}
|(t\wedge y)C_\lambda^*(s,\rho)|\leq C_d(1+|\cot\theta|)^\lambda.
\end{equation}
Combining \eqref{A-angular}, \eqref{B-angular}, and \eqref{C-angular}, we conclude that
$$|K_-(t,y)|\leq C_d(1+|\cot\theta|)^\lambda.$$
This implies \eqref{angularkernel}.
\end{proof}

\begin{remark}
The singularity in the angular estimate \eqref{angularkernel} does not correspond to an actual singularity of the Clifford-Fourier kernel, which is smooth also for parallel vectors. Rather, it reflects the loss of uniform control as the angle approaches $0$ or $\pi$. Although we do not address the optimality of the exponent $\lambda=(d-2)/2$, its value is particularly useful since $2\lambda=d-2$, which exactly compensates for the angular factor in the spherical measure. This observation is crucial for the uniform $L^2$ estimate in Lemma~\ref{est11}.
\end{remark}

\subsection{Generalized translation operator}

One of the main differences between classical Fourier analysis and Clifford-Fourier analysis lies in the definition of the translation operator. This operator was introduced in \cite{DeBieCFT} in terms of the Clifford-Fourier transform. In this part of the section we will also prove some new properties of this operator.

\begin{definition}
\label{transl}
Let $f\in\mathcal{S}(\mathbb{R}^d;\mathbb{R}_{0,d})$. For $x\in\mathbb{R}^d$, the generalized translation operator $f\mapsto\tau_x f$ is defined by
$$ \mathcal F_-(\tau_x f)(\xi) = K_-(x,\xi)\mathcal F_-f(\xi), \qquad \xi\in\mathbb R^d.$$
\end{definition}

By using the inversion formula for the CFT, see \eqref{inverse}, we can express
the generalized translation operator as
\begin{equation}
\label{tranl}(\tau_x f)(t)=(2\pi)^{-d/2}\int_{\mathbb R^d}\widetilde{K_-}(\xi,t)K_-(x,\xi)\mathcal F_-f(\xi)\,d\xi,
\qquad x,t\in\mathbb R^d.
\end{equation}
The above operator reduces to the ordinary translation operator for $f\in\mathcal{S}(\mathbb{R}^d)$ in the case $d=2$ or when $f$ is a radial function (see \cite[Theorem~7.3]{DeBieCFT}), namely
\begin{equation}
\label{radial}
\tau_x f(t)=f(t-x),\qquad t,x\in\mathbb R^d .
\end{equation}

\begin{remark}
By \cite[Proposition~3.4]{DeBieCFT}, the inverse kernel satisfies
$$\widetilde{K_-}(\xi,t)=K_-(\xi,t)$$
in even dimensions, since the Clifford-Fourier kernel has real Clifford coefficients.
Consequently, the generalized translation operator can be written as
$$(\tau_xf)(t)=(2\pi)^{-d/2}\int_{\mathbb R^d} K_-(\xi,t)K_-(x,\xi)\mathcal F_-f(\xi)\,d\xi.$$
\end{remark}

\begin{proposition}\label{mixed-kernel}
	Let $d\geq 2$ be even. For every pair of multi-indices $\alpha,\beta\in\mathbb N_0^d$,
	there exist constants $C_{\alpha,\beta}>0$ and $N_{\alpha,\beta}\geq0$ such that
	$$  \left| \partial_t^\alpha\partial_y^\beta K_-(t,y) \right| \leq C_{\alpha,\beta}
	(1+|t|+|y|)^{N_{\alpha,\beta}}, \qquad t,y\in\mathbb R^d.$$
\end{proposition}

\begin{proof}
For $d=2$, the conclusion follows directly from the explicit expression \eqref{dim2}. Assume now that $d\geq4$ is even. By the finite representation of the kernel, $K_-(t,y)$ is a finite sum of terms of the form
$$
H(t,y)=P(t,y)F_\nu(q(t,y)),
$$
where $P$ is a Clifford-valued polynomial, $\nu\geq1/2$ and we set
$$ q(t,y):=|t\wedge y|^2=|t|^2|y|^2-\langle t,y\rangle^2, \qquad F_\nu(u):=u^{-\nu/2}J_\nu(\sqrt{u}),
\qquad u>0.$$
We first show that all derivatives of $F_\nu$ are uniformly bounded on $[0,\infty)$. The function $F_\nu$ extends smoothly to $[0,\infty)$. This follows from the power series expansion of the Bessel function. Moreover, using the identity
$$ \frac{d}{dr}\bigl(r^{-\nu}J_\nu(r)\bigr)=-r^{-\nu}J_{\nu+1}(r),$$
see \cite[Formula 9.1.30]{AS}, we obtain
$$ 	F_\nu'(u)=\frac{1}{2\sqrt{u}}\frac{d}{dr}\left(r^{-\nu}J_\nu(r)\right)\bigg|_{r=\sqrt{u}}=-\frac12 u^{-(\nu+1)/2}J_{\nu+1}(\sqrt{u})=-\frac12F_{\nu+1}(u).$$
Iterating this identity gives
\begin{equation}
	\label{der1}
	F_\nu^{(k)}(u) =\frac{(-1)^k}{2^k}F_{\nu+k}(u),	\qquad k\in\mathbb N_0.
\end{equation}
The power series of $J_\nu$ ensures that these identities extend to $u=0$. The power series expansion of the Bessel function gives
$$
F_{\nu+k}(u)=\sum_{j=0}^{\infty}\frac{(-1)^j u^j}{2^{2j+\nu+k}j!\,\Gamma(j+\nu+k+1)}.
$$
Hence $F_{\nu+k}(u)$ is continuous and therefore bounded on $[0,1]$. On the other hand, the standard asymptotic estimate for Bessel functions (see \cite[Formula 9.2.5]{AS}) yields
$$|J_{\nu+k}(r)|\leq C_{\nu,k}r^{-1/2}, \qquad r\geq1.$$
Consequently, for $u\geq1$, we have
$$|F_{\nu+k}(u)|=u^{-(\nu+k)/2}|J_{\nu+k}(\sqrt{u})|\leq C_{\nu,k} u^{-(\nu+k)/2-1/4}\leq C_{\nu,k}.$$
Thus $F_{\nu+k}(u)$ is bounded on $[0,\infty)$ for every $k\geq0$. By \eqref{der1} we get
\begin{equation}\label{derbounded}
	\sup_{u\geq0}|F_\nu^{(k)}(u)|<\infty,
	\qquad k\in\mathbb N_0.
\end{equation}
We now estimate the mixed derivatives of $F_\nu(q(t,y))$. Let $\alpha,\beta\in\mathbb N_0^d$
and set $n=|\alpha|+|\beta|$. For $n=0$, the required estimate follows immediately from \eqref{derbounded}. Suppose $n\geq1$. By repeated application of the multivariate chain rule, the derivative $\partial_t^\alpha\partial_y^\beta F_\nu(q(t,y))$ is a finite linear combination of terms of the form
\begin{equation}
\label{chainF}
F_\nu^{(k)}(q(t,y))\prod_{j=1}^{k}\partial_t^{\alpha_j}\partial_y^{\beta_j}q(t,y),
\end{equation}
where $1\leq k\leq n$ and the multi-indices satisfy
$$\sum_{j=1}^{k}\alpha_j =\alpha, \qquad \sum_{j=1}^{k} \beta_j=\beta, \qquad |\alpha_j|+|\beta_j|\geq1 .$$

Since $q(t,y)$ is a homogeneous polynomial of total degree four, every mixed derivative $\partial_t^{\alpha_j}\partial_y^{\beta_j}q(t,y)$ is either identically zero or a
homogeneous polynomial of degree $4-|\alpha_j|-|\beta_j|$. Consequently, for every nonzero derivative, there exists a constant
$C_{\alpha_j,\beta_j}>0$ such that
$$ \left| \partial_t^{\alpha_j} \partial_y^{\beta_j}q(t,y) \right| \leq C_{\alpha_j,\beta_j}
(1+|t|+|y|)^{4-|\alpha_j|-|\beta_j|}.
$$
Multiplying these estimates, for each nonvanishing term we obtain
$$
\prod_{j=1}^{k} \left| \partial_t^{\alpha_j} \partial_y^{\beta_j}q(t,y)\right|
\leq C_{\alpha,\beta} (1+|t|+|y|)^{\sum_{j=1}^{k} (4-|\alpha_j|-|\beta_j|)}=C_{\alpha,\beta}
(1+|t|+|y|)^{4k-n}\leq C_{\alpha,\beta}(1+|t|+|y|)^{3n},
$$
where we have used the fact that $k\leq n$. Combining this estimate with \eqref{derbounded}, we deduce 
\begin{equation}
\label{boundcomposition}
\left| \partial_t^\alpha\partial_y^\beta F_\nu(q(t,y)) \right| \leq C_{\alpha,\beta,\nu} (1+|t|+|y|)^{3(|\alpha|+|\beta|)}.
\end{equation}
Finally, since $P(t,y)$ is a Clifford-valued polynomial, its derivatives are also polynomials. Applying Leibniz's rule to $H(t,y)=P(t,y)F_\nu(q(t,y))$ and using \eqref{boundcomposition},
we obtain constants $C_{\alpha,\beta}>0$ and $N_{\alpha,\beta}\geq0$ such that
$$ |\partial_t^\alpha\partial_y^\beta H(t,y)| \leq C_{\alpha,\beta} (1+|t|+|y|)^{N_{\alpha,\beta}}.$$
Since $K_-(t,y)$ is a finite sum of terms of this form, the desired estimate follows.
\end{proof}

\begin{proposition}
\label{tranS}
Let $d\geq2$ be even and let $x\in\mathbb R^d$ be fixed. Then the generalized translation preserves the Schwartz space.
\end{proposition}
\begin{proof}
By Proposition \ref{mixed-kernel} we deduce that for every multi-index $\beta$, there exist constants $C_{\beta,x}>0$ and $N_\beta\geq0$ such that
$$ \bigl|\partial_\xi^\beta K_-(x,\xi)\bigr|
\leq C_{\beta,x}(1+|\xi|)^{N_\beta}.$$
By Leibniz's rule, for every multi-index $\alpha$, we have
$$ \left| \xi^\alpha \partial_\xi^\beta \bigl( K_-(x,\xi)\mathcal F_-(f)(\xi)\bigr)\right| \leq C_{\alpha,\beta,x} \sum_{\gamma\leq\beta} (1+|\xi|)^{|\alpha|+N_\gamma} 	
\left| \partial_\xi^{\beta-\gamma}\mathcal F_-(f)(\xi)\right|.
$$
The right-hand side is bounded since $\mathcal F_-(f)$ is a Schwartz
function. Therefore,
$$
K_-(x,\cdot)\mathcal F_-(f) \in
\mathcal S(\mathbb R^d;\mathbb R_{0,d}).
$$
Hence, by Definition~\ref{transl} and the fact that the Clifford-Fourier
transform maps the Schwartz space into itself in even dimension, we obtain $ \tau_x f \in \mathcal S(\mathbb R^d;\mathbb R_{0,d})$.
\end{proof}

In the sequel, we will need a uniform $L^2$-estimate for the generalized translation operator. To establish this estimate, we first prove the following result.

\begin{lemma}
\label{translationest}
Let $d\geq 2$ be even and let $g\in\mathcal S(\mathbb R^d;\mathbb R_{0,d})$. Then, for every $t\in\mathbb R^d$, there exists a constant $C_d>0$, depending only on the dimension, such that
$$ \int_{\mathbb R^d}|(\tau_xg)(t)|^2\,dx \leq C_d \int_{\mathbb R^d} |K_-(\xi,t)|^2|\mathcal F_-g(\xi)|^2\,d\xi.$$
Moreover, we have
$$ \int_{\mathbb R^d}|(\tau_xg)(t)|^2\,dx \leq C_d(1+|t|)^{2\lambda}\mathcal A_g,$$
where
\begin{equation}
	\label{new}
	\mathcal A_g:= \int_{\mathbb R^d} (1+|\xi|)^{2\lambda} |\mathcal F_-(g)(\xi)|^2\,d\xi.
\end{equation}
\end{lemma}

\begin{proof}
Since both the Clifford-Fourier kernel and the Clifford-Fourier transform
take values in the Clifford algebra $\mathbb R_{0,d}$, see Definition~\ref{cft} and formula  \eqref{kernelreal}, respectively, we can write
$$ K_-(\xi,t)=\sum_A a_A(\xi,t)e_A, \qquad
\mathcal F_-g(\xi)=\sum_B b_B(\xi)e_B,$$
where $a_A(\xi,t),b_B(\xi)\in\mathbb R$. Thus, by \eqref{tranl}, we can express the generalized translation operator as
$$ (\tau_xg)(t) = \sum_{A,B} e_A \left[ \frac{1}{(2\pi)^{d/2}} \int_{\mathbb R^d} K_-(x,\xi)a_A(\xi,t)b_B(\xi)\,d\xi \right] e_B.$$
For a real-valued scalar function $h$,  we set
$$ Th(x):= \frac{1}{(2\pi)^{d/2}}\int_{\mathbb R^d}K_-(x,\xi)h(\xi)\,d\xi.$$
By \cite[Prop. 3.5]{LiLengFeiPW} we have $K_-(x,\xi)=\overline{K_-(\xi,x)}$. This implies that
$$ Th(x)=\overline{\mathcal F_-h(x)}.$$
Therefore, by Theorem \ref{par}, we get
$$ \|Th\|_{L^2(\mathbb R^d;\mathbb R_{0,d})} =\|h\|_{L^2(\mathbb R^d)}.$$
By \eqref{estmou} there exists a constant $C_d>0$ such that
\begin{eqnarray*}
\nonumber
\int_{\mathbb R^d}|(\tau_xg)(t)|^2\,dx &\leq& C_d \sum_{A,B} \int_{\mathbb R^d} |T(a_Ab_B)(x)|^2\,dx= C_d \sum_{A,B} \int_{\mathbb R^d} |a_A(\xi,t)|^2|b_B(\xi)|^2\,d\xi\\
\label{N1}
&=& C_d \int_{\mathbb R^d} |K_-(\xi,t)|^2|\mathcal F_-g(\xi)|^2\,d\xi.
\end{eqnarray*}
Finally, by Proposition \ref{est} we get
$$ 	\int_{\mathbb R^d}|(\tau_xg)(t)|^2\,dx \leq C_d(1+|t|)^{2\lambda} \int_{\mathbb R^d}
(1+|\xi|)^{2\lambda} |\mathcal F_-g(\xi)|^2\,d\xi=
C_d(1+|t|)^{2\lambda}\mathcal A_g.$$
\end{proof}

\begin{lemma}
	\label{est11}
	Let $d\geq 2$ be even and let $h\in\mathcal S(\mathbb R^d;\mathbb R_{0,d})$. Then there exists a constant $C_h>0$ such that
	$$\int_{\mathbb R^d}|(\tau_xh)(t)|^2\,dx \leq C_h$$
	for every $t\in\mathbb R^d$.
\end{lemma}

\begin{proof}
	For $d=2$, the generalized translation coincides with the ordinary translation, see \eqref{radial}. Hence the conclusion follows immediately. Assume now that $d\geq4$ is even. Since $\mathcal F_-h$ belongs to the Schwartz space, for every $L>0$ there exists $M_{h,L}>0$
	such that
	$$ |\mathcal F_-h(\xi)| \leq M_{h,L}(1+|\xi|)^{-L}. $$
	Choose $L>d/2$. For $t\neq0$, by Lemma \ref{translationest} and passing to polar coordinates $\xi=\rho\eta$, with $\rho>0$ and $\eta\in \mathbb{S}^{d-1}$, we obtain
	\begin{equation}
		\label{esttra}
		\int_{\mathbb R^d}|(\tau_xh)(t)|^2\,dx \leq C_dM_{h,L}^2
		\int_0^\infty \frac{\rho^{d-1}}{(1+\rho)^{2L}}\,d\rho \int_{\mathbb{S}^{d-1}} \frac{1}{(\sin\theta)^{2\lambda}}d\sigma(\eta).
	\end{equation}
	Since $2\lambda=d-2$ and $d\sigma(\eta)=(\sin\theta)^{d-2}\,d\theta\,d\sigma_{\mathbb{S}^{d-2}}$ we have
	$$ \int_{\mathbb{S}^{d-1}} \frac{1}{(\sin\theta)^{2\lambda}} d\sigma(\eta)=\pi|\mathbb{S}^{d-2}|<\infty.
	$$
	Moreover, since $2L>d$, we have
	$$ \int_0^\infty \frac{\rho^{d-1}}{(1+\rho)^{2L}}\,d\rho <\infty.$$
	By \eqref{esttra} we have
	$$ \int_{\mathbb R^d}|(\tau_xh)(t)|^2\,dx \leq C_h, $$
	where $C_h$ is independent of $t\neq0$. The case $t=0$ follows from $K_-(\xi,0)=1$, Lemma \ref{translationest} and Theorem \ref{par}. Thus the conclusion holds for every $t\in\mathbb R^d$.
\end{proof}

\section{The Clifford short-time Fourier transform with non-radial windows}

The short-time Fourier transform is a powerful tool in time-frequency analysis for obtaining information about the local properties of a signal. The main idea is to localize the signal around a given point and then compute the Fourier transform of the localized signal. This localization is usually achieved by multiplying the signal by a smooth function, called the window function. The theory of the short-time Fourier transform is very rich and has been systematically developed in \cite{G}. Connections with polyanalytic Bargmann-Fock spaces and superoscillations can be found in \cite{AB, ADDS}.

\medskip

The notion of the short-time Fourier transform has also been generalized to the hypercomplex setting in several papers; see, for example, \cite{DK, DK1}. In this paper, we follow the approach introduced in \cite{DeMartinoCSTFT}, where the author defined a Clifford short-time Fourier transform based on the Clifford-Fourier transform recalled in Section~3. However, unlike in \cite{DeMartinoCSTFT}, we do not restrict ourselves to radial window functions. Our main goal is to investigate the structure of this transform and establish some uncertainty principles associated with it.

\begin{definition}
\label{CSFT}
Let $f$, $g \in \mathcal{S}(\mathbb{R}^d; \mathbb{R}_{0,d})$. The Clifford short-time Fourier transform (CSTFT) of a function $f$ with respect to a \emph{window function} $g$ is defined as
$$V_g f(x, \omega)= \mathcal{F}_{-}(\tau_x g \cdot f)(\omega), \qquad x, \omega \in \mathbb{R}^d,$$
where $ \mathcal{F}_{-}$ is the Clifford Fourier transform, see Definition \ref{cft}.
\end{definition}

\begin{remark}
If the function $g$ is radial we can write the CSTFT as
\begin{equation}
	\label{St1}
	V_g f(x,\omega)=(2 \pi)^{-\frac{d}{2}}\int_{\mathbb{R}^d}K_{-}(t,\omega)g(t-x)f(t)\,dt ,\qquad x,\omega\in\mathbb R^d.
\end{equation}
This way of writing the CSTFT follows by \eqref{radial}. However if the function $g$ is not radial we can write the CSTFT as
\begin{equation}
	\label{St2}
	V_g f(x,\omega)=(2 \pi)^{-\frac{3}{2}d} \int_{\mathbb{R}^{3d}}K_{-}(t, \omega)K_{-}(\xi,t)K_{-}(x,\xi)K_{-}(s, \xi) g(s) f(t)dsd \xi dt.
\end{equation}
\end{remark}

In this section, we aim to investigate the main properties of the CSTFT when the window function is not radial.

\begin{theorem}\label{wellD}
Let $d>2$ be even and let $f\in L^2(\mathbb R^d;\mathbb R_{0,d})$ and $g\in\mathcal S(\mathbb R^d;\mathbb R_{0,d})$. Then $V_gf(x,\omega)$ is well defined for every $x,\omega\in\mathbb R^d$.
\end{theorem}

\begin{proof}
By Definition~\ref{CSFT}, the CSTFT is the Clifford-Fourier transform of $(\tau_xg)\cdot f$. Thus, by Theorem~\ref{cont1}, it suffices to prove that $(\tau_xg)\cdot f\in\mathcal B(\mathbb R^d;\mathbb R_{0,d})$. For every fixed $x\in\mathbb R^d$, by Proposition~\ref{tranS} we have that $\tau_xg\in\mathcal S(\mathbb R^d;\mathbb R_{0,d})$. In particular, setting $\lambda=(d-2)/2$,
we have
$$ (1+|\cdot|)^\lambda\tau_xg \in L^2(\mathbb R^d;\mathbb R_{0,d}).$$
By H\"older's inequality and \eqref{estmou}, we obtain
$$
\begin{aligned}
\|(\tau_xg)f\|_{\mathcal B} &=\int_{\mathbb R^d}(1+|t|)^\lambda |(\tau_xg)(t)f(t)|\,dt\leq C_d
\int_{\mathbb R^d}(1+|t|)^\lambda |(\tau_xg)(t)|\,|f(t)|\,dt\\
&\leq C_d \left(\int_{\mathbb R^d}(1+|t|)^{2\lambda} |(\tau_xg)(t)|^2\,dt\right)^{1/2}\|f\|_2<\infty.
\end{aligned}
$$
	Therefore, $(\tau_xg)f$ belongs to
	$\mathcal B(\mathbb R^d;\mathbb R_{0,d})$,
	which completes the proof.
\end{proof}

Now, we study the main properties of the CSTFT. In order to do this we need the following notation
\begin{equation}
\label{Wintegral}
W_{g_1,g_2}(t):=\int_{\mathbb R^d}\overline{(\tau_xg_1)(t)} (\tau_xg_2)(t)\,dx, \quad t \in \mathbb{R}^d, \quad g_1, g_2 \in \mathcal S(\mathbb R^d; \mathbb{R}_{0,d}).
\end{equation}

\begin{theorem}[Weighted Moyal identity]
\label{wMoyal}
Let $d>2$ be even and $\lambda=\frac{d-2}{2}$. Let $g_1$, $g_2 \in\mathcal S(\mathbb R^d; \mathbb{R}_{0,d})$, $f_1$, $f_2 \in L^2(\mathbb{R}^d; \mathbb{R}_{0,d})$. Then we have
\begin{equation}
\label{WMoyal1}
\int_{\mathbb R^{2d}} \overline{V_{g_1}f_1(x,\omega)} V_{g_2}f_2(x,\omega)\,d\omega\,dx=
\int_{\mathbb R^d}\overline{f_1(t)} W_{g_1,g_2}(t) f_2(t)\,dt.
\end{equation}
Moreover, both integrals are absolutely convergent.
\end{theorem}

\begin{proof}
For every $x\in\mathbb R^d$, by  Theorem \ref{par}, we have
$$
\int_{\mathbb R^d}\overline{V_{g_1}f_1(x,\omega)}V_{g_2}f_2(x,\omega)\,d\omega=\int_{\mathbb R^d}\overline{(\tau_xg_1)(t)f_1(t)}(\tau_xg_2)(t)f_2(t)\,dt.
$$
By \eqref{reverse} we have
$$\int_{\mathbb R^{2d}} \overline{V_{g_1}f_1(x,\omega)}V_{g_2}f_2(x,\omega)\,d\omega\,dx=
\int_{\mathbb R^d}\overline{f_1(t)} \left(\int_{\mathbb R^d}\overline{(\tau_xg_1)(t)} (\tau_xg_2)(t)\,dx\right)f_2(t)\,dt,
$$
which proves \eqref{WMoyal1}. It remains to prove that the two integrals in \eqref{WMoyal1} are absolutely convergent. By \eqref{estmou}, Theorem \ref{par}, Tonelli's theorem and Lemma \ref{est11} there exists $C_{g_j}>0$ such that
\begin{eqnarray*}
\int_{\mathbb R^{2d}} |V_{g_j}f_j(x,\omega)|^2\,d\omega\,dx&=& \int_{\mathbb R^{2d}} |(\tau_xg_j)(t)f_j(t)|^2\,dtdx \leq 	C_d\int_{\mathbb R^d}|f_j(t)|^2 \left( \int_{\mathbb R^d}|(\tau_xg_j)(t)|^2\,dx \right)dt\\
&\leq &C_d C_{g_j}  \int_{\mathbb R^d} |f_j(t)|^2\,dt=
C_d C_{g_j} \|f_j\|_{L^2}^2 <\infty.
\end{eqnarray*}
This implies that $V_{g_j}f_j\in L^2(\mathbb R^{2d};\mathbb R_{0,d})$, for $j=1,2$. Thus, by the Cauchy-Schwarz inequality, the integral on the left-hand side of \eqref{WMoyal1} is absolute convergent. We now consider the right-hand side. By the definition of
$W_{g_1,g_2}$, \eqref{estmou}, and the Cauchy-Schwarz inequality in the $x$ variable, and Lemma \ref{est11} there exist $C_{g_1}$, $C_{g_2}>0$ such that
$$ |W_{g_1,g_2}(t)| \leq C_d \left( \int_{\mathbb R^d}|(\tau_xg_1)(t)|^2\,dx
\right)^{1/2} \left( \int_{\mathbb R^d}|(\tau_xg_2)(t)|^2\,dx \right)^{1/2} \leq  C_{d} C_{g_1}^{\frac{1}{2}} C_{g_2}^{\frac{1}{2}}.$$
 Thus we have
$$ \int_{\mathbb R^d} \left| \overline{f_1(t)}W_{g_1,g_2}(t)f_2(t) \right|\,dt \leq
C_{d} C_{g_1}^{\frac{1}{2}} C_{g_2}^{\frac{1}{2}} \int_{\mathbb R^d} |f_1(t)|\,|f_2(t)|\,dt \leq C_{d} C_{g_1}^{\frac{1}{2}} C_{g_2}^{\frac{1}{2}} 	\|f_1\|_{L^2} \|f_2\|_{L^2} <\infty,
$$
where in the last inequality we used Cauchy-Schwarz. Hence also the integral on the right-hand side of \eqref{WMoyal1} is absolutely convergent.
\end{proof}

\begin{remark}
If $g_1$ and $g_2$ are real-valued radial windows, then
$$ W_{g_1,g_2}(t)=\int_{\mathbb R^d}g_1(z)g_2(z)\,dz,$$
which is independent of $t$. Therefore, Theorem~\ref{wMoyal} reduces to the usual orthogonality relation for the CSTFT; see~\cite[Thm.~8.1]{DeMartinoCSTFT}.
\end{remark}

\begin{remark}
\label{specificcase}
In the case $g_1=g_2=g$ in \eqref{Wintegral} we denote the integral by $W_g(t)$.
\end{remark}
\begin{proposition}\label{SchwartzCSTFT}
	Let $d\geq2$ be even and let $f$, $g\in\mathcal S(\mathbb R^d;\mathbb R_{0,d})$. Then we have
	$V_gf\in\mathcal S(\mathbb R^{2d};\mathbb R_{0,d})$.
\end{proposition}

\begin{proof}
We start writing 
$$ K_-(\xi,t)=\sum_A a_A(\xi,t)e_A, \qquad \mathcal F_-g(\xi)=\sum_B b_B(\xi)e_B, \qquad f(t)=\sum_C f_C(t)e_C, $$
where $b_B(\xi)$, and $f_C(t)$ are real-valued scalar functions. Since $f\in\mathcal S(\mathbb R^d;\mathbb R_{0,d})$ and $\mathcal F_-g\in\mathcal S(\mathbb R^d;\mathbb R_{0,d})$, we have $f_C,b_B\in\mathcal S(\mathbb R^d)$. By Leibniz's rule, every mixed derivative $
 \partial_\xi^\alpha\partial_t^\beta \bigl(a_A(\xi,t)b_B(\xi)f_C(t)\bigr)$
is a finite sum of terms of the form
$$ \bigl(\partial_\xi^{\alpha_1} \partial_t^{\beta_1}a_A(\xi,t)\bigr) \bigl(\partial_\xi^{\alpha_2}b_B(\xi)\bigr) \bigl(\partial_t^{\beta_2}f_C(t)\bigr),
$$
where $\alpha_1+\alpha_2=\alpha$ and $\beta_1+\beta_2=\beta$. By Proposition~\ref{mixed-kernel}, there exist constants $C>0$ and $N\geq0$ such that
$$ \left| \partial_\xi^{\alpha_1} \partial_t^{\beta_1}a_A(\xi,t) \right| \leq C(1+|\xi|+|t|)^N.
$$
On the other hand, since $b_B$ and $f_C$ are Schwartz functions, for every $M\geq0$
there exists a constant $C_M>0$ such that
$$ \left|\partial_\xi^{\alpha_2}b_B(\xi)\right| \leq C_M(1+|\xi|)^{-M},
\qquad \left|\partial_t^{\beta_2}f_C(t)\right| \leq C_M(1+|t|)^{-M}.
$$
Consequently,
$$\left|\xi^\mu t^\nu \bigl(\partial_\xi^{\alpha_1} \partial_t^{\beta_1}a_A(\xi,t)\bigr)
\bigl(\partial_\xi^{\alpha_2}b_B(\xi)\bigr) \bigl(\partial_t^{\beta_2}f_C(t)\bigr)\right|\leq
CC_M^2 (1+|\xi|+|t|)^N (1+|\xi|)^{|\mu|-M} (1+|t|)^{|\nu|-M}$$
Using the fact that $ 1+|\xi|+|t| \leq(1+|\xi|)(1+|t|)$, and choosing $M\geq N+|\mu|+|\nu|$,
the right-hand side is uniformly bounded in $(\xi,t)$. Since the Leibniz rule involves only finitely many terms, we conclude that
$$ \sup_{\xi,t\in\mathbb R^d} \left| \xi^\mu t^\nu
\partial_\xi^\alpha\partial_t^\beta \bigl(a_A(\xi,t)b_B(\xi)f_C(t)\bigr) \right|<\infty .$$

Thus, setting $h_{A,B,C}(\xi,t):=a_A(\xi,t)b_B(\xi)f_C(t)$, using the identity
$K_-(x,\xi)=\overline{K_-(\xi,x)}$, and by \eqref{tranl} we obtain
\begin{equation}\label{star}
(\tau_xg)(t)f(t)=\sum_{A,B,C}e_A\,\overline{(\mathcal F_-^\xi h_{A,B,C})(x,t)}\,e_Be_C.
\end{equation}
Here $\mathcal F_-^\xi$ denotes the Clifford-Fourier transform with respect to the variable $\xi$. We now show that the partial Clifford-Fourier transform preserves the Schwartz space on $\mathbb R^{2d}$. Let $h\in\mathcal S(\mathbb R^{2d})$. Since the kernel of the transform is independent of $t$, differentiation with respect to $t$ commutes with the partial transform. Thus, for every pair of multi-indices $\beta,\nu$, we have
\begin{equation}
\label{change}
t^\nu\partial_t^\beta \mathcal F_-^\xi h(x,t)= \mathcal F_-^\xi
\bigl(t^\nu\partial_t^\beta h\bigr)(x,t).
\end{equation}

By the continuity of the Clifford-Fourier transform on the Schwartz space
(Theorem~\ref{cont1}), for every pair of multi-indices $\alpha,\mu$ there exist
a constant $C>0$ and finitely many pairs of multi-indices $(\gamma_j,\delta_j)$
such that
$$\sup_{x\in\mathbb R^d}
\left| x^\mu\partial_x^\alpha \mathcal F_-v(x) \right| \leq
C\sum_{j=1}^{M} \sup_{\xi\in\mathbb R^d}
\left| \xi^{\gamma_j} \partial_\xi^{\delta_j}v(\xi) \right|, \qquad \forall v\in\mathcal S(\mathbb R^d)
$$
Applying the above estimate to $v(\xi)=t^\nu\partial_t^\beta h(\xi,t)$, with $t$ fixed, using \eqref{change} and then taking the supremum with respect to $t$, we obtain
$$ \sup_{x,t\in\mathbb R^d}
\left| x^\mu t^\nu \partial_x^\alpha\partial_t^\beta (\mathcal{F}_{-}^\xi h)(x,t) \right| \leq
C\sum_{j=1}^{M} \sup_{\xi,t\in\mathbb R^d}
\left| \xi^{\gamma_j}t^\nu \partial_\xi^{\delta_j}\partial_t^\beta h(\xi,t)
\right|<\infty.$$
The last inequality follows from the fact that $h\in\mathcal S(\mathbb R^{2d})$. Therefore, $\mathcal F_-^\xi h \in\mathcal S(\mathbb R^{2d};\mathbb R_{0,d})$, and so we have
$ \overline{\mathcal F_-^\xi h_{A,B,C}}\in\mathcal S(\mathbb R^{2d};\mathbb R_{0,d}).$ Therefore, by \eqref{star} we have $(\tau_xg)(t)f(t) \in \mathcal S(\mathbb R^{2d};\mathbb R_{0,d})$. The conclusion follows from Definition~\ref{CSFT} and the fact that the Clifford-Fourier transform preserves the Schwartz space.

\end{proof}

We now introduce a class of Clifford-valued functions that will be suitable for establishing a reconstruction formula for the CSTFT with a non-radial window function.

\begin{definition}
\label{adm}
Let $g\in\mathcal S(\mathbb R^d;\mathbb R_{0,d})$.  We set

\begin{equation}
	\label{Qinte}
	Q_g(t,u):= \int_{\mathbb{R}^{d}} | (\tau_x g)(t) u|^2 dx, \qquad
	t\in\mathbb R^d,\quad u\in\mathbb R_{0,d}.
\end{equation}
We say that $g$ is an admissible window if there exist constants $c_g>0$ and $N\geq0$ such that
\begin{equation}
\label{inedef}
Q_g(t,u) \geq c_g(1+|t|)^{-N}|u|^2
\end{equation}
for every $t\in\mathbb R^d$ and every $u\in\mathbb R_{0,d}$.
\end{definition}

\begin{remark}
Every nonzero real-valued radial function
$g\in\mathcal S(\mathbb R^d)$ is admissible. In fact, by \eqref{radial} we have, for every $u\in\mathbb R_{0,d}$, that
$$ Q_g(t,u)=\|g\|_2^2|u|^2.$$
Hence the admissibility condition holds with $N=0$ and $c_g=\|g\|_2^2$.
\end{remark}

We now show that admissibility is stable under sufficiently small perturbations of nonzero real-valued radial windows.

\begin{proposition}
	\label{exx}
	Let $d\geq2$ be even. Let $\varphi\in\mathcal S(\mathbb R^d)$ be a nonzero real-valued radial function and let $h\in\mathcal S(\mathbb R^d;\mathbb R_{0,d})$. Then there exists $\varepsilon_0>0$ such that
	$$ g_\varepsilon=\varphi+\varepsilon h$$
	is an admissible window, with $N=0$, for every $|\varepsilon|<\varepsilon_0$.
\end{proposition}

\begin{proof}
	Since $\varphi$ is radial, by \eqref{radial} we have $(\tau_x\varphi)(t)=\varphi(t-x)$. By linearity, we have
	$$(\tau_xg_\varepsilon)(t)=\varphi(t-x)+\varepsilon(\tau_xh)(t).$$
	By the reverse triangle inequality in $L^2$ and Lemma \ref{est11} we have
	$$\left(\int_{\mathbb R^d}|(\tau_xg_\varepsilon)(t)u|^2\,dx\right)^{1/2}\geq
	\|\varphi\|_2|u|-|\varepsilon|\left(\int_{\mathbb R^d}|(\tau_xh)(t)u|^2\,dx\right)^{1/2}
	\geq \bigl(\|\varphi\|_2-|\varepsilon|D_h\bigr)|u|,$$
	where $D_h>0$ is independent of $t$ and $u$. Choose
	$$ \varepsilon_0=\frac{\|\varphi\|_2}{2D_h}. $$
	For $|\varepsilon|<\varepsilon_0$, we get
	$$ \int_{\mathbb R^d} |(\tau_xg_\varepsilon)(t)u|^2\,dx \geq\frac{\|\varphi\|_2^2}{4}|u|^2.$$
	Therefore, $g_\varepsilon$ is admissible with $N=0$.
\end{proof}

Now, we provide a specific example of a non-radial admissible function.

\begin{example}\label{non-radial-gaussian}
	Let $d\geq2$ be even and consider
	$$ g_\varepsilon(t) =e^{-|t|^2/2}(1+\varepsilon t_1e_1),\qquad t\in\mathbb R^d.$$
	The function $g_\varepsilon$ is non-radial for every $0<|\varepsilon|<\varepsilon_0$. We set
	$$ \varphi(t):=e^{-|t|^2/2}, \qquad h(t):=t_1e^{-|t|^2/2}e_1.$$
	Clearly, $\varphi$ is a nonzero real-valued radial Schwartz function and $h\in\mathcal S(\mathbb R^d;\mathbb R_{0,d})$. Thus, by Proposition~\ref{exx} the function $g_\varepsilon$ is admissible with $N=0$.
\end{example}

We now establish an upper bound for the quadratic expression defined in \eqref{Qinte}.

\begin{proposition}
\label{Qgestimate}
Let $d\geq2$ be even and let $g\in\mathcal S(\mathbb R^d;\mathbb R_{0,d})$. If $d=2$, then
\begin{equation}
\label{d2}
Q_g(t,u)=\|g\|_{L^2(\mathbb R^2; \mathbb{R}_{0,2})}^2\,|u|^2, \qquad t\in\mathbb R^2,\quad u\in\mathbb R_{0,2}.
\end{equation}
If $d\geq4$, we have
\begin{equation} \label{Qgtwosided}
Q_g(t,u) \leq c_d \mathcal{A}_g(1+|t|)^{d-2}|u|^2,
\end{equation}
where $\mathcal{A}_g$ has been introduced in \eqref{new} and $c_d$ is a constant depending on the dimension.
\end{proposition}

\begin{proof}
Assume first that $d=2$. In this case the generalized translation coincides with the classical translation: $(\tau_xg)(t)=g(t-x)$. Moreover, the norm on $\mathbb R_{0,2}$ is multiplicative. Hence

$$ Q_g(t,u) = \int_{\mathbb R^2}|g(t-x)|^2\,dx\,|u|^2 = \|g\|_{L^2(\mathbb R^2; \mathbb{R}_{0,2})}^2|u|^2.$$
Let now $d\geq4$.  By Lemma \ref{translationest} we get
$$Q_g(t,u) \leq 2^d |u|^2 (1+|t|)^{2\lambda} \mathcal{A}_g.$$
 Since we have $2\lambda=d-2$ we get the upper estimate of \eqref{Qgtwosided}.
\end{proof}

\begin{lemma}
\label{adm1}
Let $g\in\mathcal S(\mathbb R^d;\mathbb R_{0,d})$ be an admissible window. Then $W_g(t)$ is invertible for every $t\in\mathbb R^d$ and there exist $c_g>0$ and $N \geq 0$ such that
\begin{equation}
\label{inequp}
|W_g(t)^{-1}| \leq c_g^{-1}(1+|t|)^N.
\end{equation}
\end{lemma}

\begin{proof}
By \eqref{inner1}, for every $u\in\mathbb R_{0,d}$, we have
\begin{equation}
\label{inner}
(u,W_g(t)u)= \operatorname{Sc}\left( \overline u \int_{\mathbb R^d}\overline{\tau_xg(t)}\tau_xg(t)\,dx\,u \right)=\int_{\mathbb R^d}\operatorname{Sc}\left(\overline{\tau_xg(t)u}\,\tau_xg(t)u\right)\,dx=\int_{\mathbb R^d}|\tau_xg (t)u|^2\,dx.
\end{equation}
Since by hypothesis $g$ is an admissible window, there exist constants $c_g>0$ and $N \geq 0$, see Definition \ref{adm}, such that
\begin{equation}
\label{inter}
(u,W_g(t)u) \geq c_g(1+|t|)^{-N}|u|^2.
\end{equation}
Therefore the left multiplication operator
$$ L_{W_g(t)}u=W_g(t)u$$
is injective. Since $\mathbb R_{0,d}$ is finite-dimensional, $L_{W_g(t)}$ is invertible, and consequently $W_g(t)$ is invertible.
Now, we prove \eqref{inequp}. Let $v=W_g(t)u$. By the Cauchy-Schwarz inequality and \eqref{inter} we have
$$ c_g(1+|t|)^{-N}|u|^2 \leq (u,v) \leq |u|\,|v|.$$
This implies that
$$ |u| \leq c_g^{-1}(1+|t|)^N|v|,$$
and hence by the definition of norm of an operator we have
$$ \|L_{W_g(t)}^{-1}\|_{\mathrm{op}} \leq c_g^{-1}(1+|t|)^N.$$
Finally, since $ L_{W_g(t)}^{-1}=L_{W_g(t)^{-1}}$, we obtain
$$|W_g(t)^{-1}|=|L_{W_g(t)^{-1}}1| \leq c_g^{-1}(1+|t|)^N.$$
\end{proof}

\begin{remark}
	Although the quantities $W_g$ and $Q_g$ are closely related, they play
	different roles. The function $W_g(t)$ is Clifford-valued and arises naturally
	in the weighted Moyal identity and in the reconstruction formula, whereas
	$Q_g(t,u)$ is a nonnegative real-valued quadratic form. More precisely, by \eqref{inner} we have
	$$ (u,W_g(t)u)=Q_g(t,u).$$
\end{remark}

\begin{theorem}[Reconstruction formula]
	Let $d>2$ be even and let $f$, $g\in\mathcal S(\mathbb R^d;\mathbb R_{0,d})$. We assume that $g$ is an admissible function. Then for every $y \in \mathbb{R}^d$ we have
	\begin{equation}
		\label{recon}
		f(y)=\frac{1}{(2\pi)^{d/2}} W_g(y)^{-1} \int_{\mathbb R^{2d}} \overline{(\tau_x g)(y)} K_-(\omega,y) V_gf(x,\omega) \,d\omega\,dx.
	\end{equation}
\end{theorem}

\begin{proof}
By Proposition~\ref{SchwartzCSTFT}, we have $(\tau_x g)(\cdot)f(\cdot)\in \mathcal S(\mathbb R^d;\mathbb R_{0,d})$ for every $x\in\mathbb R^d$, and $V_gf\in\mathcal S(\mathbb R^{2d};\mathbb R_{0,d})$. In particular, both $(\tau_x g)f$ and $V_gf(x,\cdot)$ belong to $\mathcal B(\mathbb R^d;\mathbb R_{0,d})$. Hence, by \cite[Thm.~8.4]{DeBieCFT}, we obtain
\begin{equation}\label{oneR}
(\tau_x g)(y)f(y)=\frac{1}{(2\pi)^{d/2}}\int_{\mathbb R^d}K_-(\omega,y)V_gf(x,\omega)\,d\omega.
\end{equation}
Multiplying \eqref{oneR} from the left by $\overline{(\tau_x g)(y)}$ gives
$$\overline{(\tau_x g)(y)}(\tau_x g)(y)f(y)=\frac{1}{(2\pi)^{d/2}} \int_{\mathbb R^d}
\overline{(\tau_x g)(y)} K_-(\omega,y)V_gf(x,\omega)\,d\omega.$$
We now justify the absolute convergence of the integral over $\mathbb R^{2d}$. Since $V_gf$ belongs to the Schwartz space, for every $M>0$ there exists a constant $C_M>0$ such that
$$ |V_gf(x,\omega)| \leq C_M(1+|x|)^{-M}(1+|\omega|)^{-M}.$$
By \eqref{estmou} and Proposition \ref{est} we have
$$ \int_{\mathbb R^{2d}} \left| \overline{(\tau_x g)(y)} K_-(\omega,y)V_gf(x,\omega) \right|\,d\omega\,dx \leq C_{d,y,M} \left(\int_{\mathbb R^d}
|(\tau_x g)(y)|(1+|x|)^{-M}\,dx\right) \left(\int_{\mathbb R^d} (1+|\omega|)^{\lambda-M}\,d\omega\right).
$$
Choose $M>d+\lambda$. Then the integral with respect to $\omega$ is finite. On the other hand, by the Cauchy-Schwarz inequality and Lemma~\ref{est11} we get
$$\int_{\mathbb R^d} |(\tau_x g)(y)|(1+|x|)^{-M}\,dx \leq \left(\int_{\mathbb R^d} |(\tau_x g)(y)|^2\,dx\right)^{1/2} \left(\int_{\mathbb R^d} (1+|x|)^{-2M}\,dx\right)^{1/2} <\infty.
$$
Thus the integral over $\mathbb R^{2d}$ is absolutely convergent. We integrate \eqref{oneR}, after multiplication from the left by $\overline{(\tau_x g)(y)}$, with respect to $x$ and using Fubini's theorem (justified from the previous arguments), we obtain
$$\left( \int_{\mathbb R^d} \overline{(\tau_x g)(y)}(\tau_x g)(y)\,dx \right)f(y)
=\frac{1}{(2\pi)^{d/2}} \int_{\mathbb R^{2d}} \overline{(\tau_x g)(y)} K_-(\omega,y)V_gf(x,\omega)\,d\omega\,dx.$$
By the definition of $W_g$, this yields
$$ W_g(y)f(y) =\frac{1}{(2\pi)^{d/2}} \int_{\mathbb R^{2d}} \overline{(\tau_x g)(y)} K_-(\omega,y)V_gf(x,\omega)\,d\omega\,dx.$$
Finally, since $g$ is admissible, Lemma~\ref{adm1} guarantees that $W_g(y)$ is invertible. Multiplying from the left by $W_g(y)^{-1}$, we obtain \eqref{recon}.
\end{proof}

\begin{remark}
The admissibility condition in Definition~\ref{adm} ensures that $W_g(t)$ is invertible for every $t\in\mathbb R^d$ and that its inverse has at most polynomial growth. This allows us to recover $f$ from its CSTFT by multiplying by $W_g(t)^{-1}$. Without admissibility, $W_g(t)$ may fail to be invertible, preventing the reconstruction procedure developed here. Thus, admissibility provides a sufficient quantitative nondegeneracy condition, although it is not claimed to be necessary for reconstruction.
\end{remark}

Now, we introduce a counterpart of the Gabor space for the CSTFT with non-radial window function. This is defined as
 $$\mathcal G_g:=\left\{ V_gf: f\in  \mathcal{S}(\mathbb R^d;\mathbb R_{0,d})
 \right\},$$

\begin{theorem}
Let $d>2$ be even and let $f$, $g\in\mathcal S(\mathbb R^d;\mathbb R_{0,d})$ and $g$ be an admissible function. Then the range $\mathcal G_g$ satisfies the following reproducing formula
\begin{equation}
	\label{repS}
	V_gf(x',\omega') = \int_{\mathbb R^{2d}} \mathcal K_g(\omega,x;\omega',x') V_gf(x,\omega) \,d\omega\,dx.
\end{equation}
where
\begin{equation}
\label{kenrel1}
\mathcal K_g(\omega,x;\omega',x')=	\frac{1}{(2\pi)^{d}}
\int_{\mathbb R^d} K_-(t,\omega') (\tau_{x'} g)(t) W_g(t)^{-1} \overline{(\tau_x g)(t)} K_-(\omega,t)\,dt.
\end{equation}
Furthermore, the integral in \eqref{kenrel1} is absolutely convergent.
\end{theorem}

\begin{proof}
By definition of the Clifford short-time Fourier transform, see Definition \ref{CSFT} and the reconstruction formula we obtain
$$ V_gf(x',\omega')= \frac{1}{(2\pi)^{d}} \int_{\mathbb R^d} K_-(t,\omega') (\tau_{x'} g)(t) W_g(t)^{-1}\left( \int_{\mathbb R^{2d}}
\overline{(\tau_x g)(t)} K_-(\omega,t) V_gf(x,\omega) \,d\omega\,dx
\right)dt
.
$$
We now justify interchanging the order of integration in the above expression. By Proposition \ref{SchwartzCSTFT}, for every $s\geq0$, we deduce that
\begin{equation}
\label{conV}
\int_{\mathbb R^{2d}}
(1+|\omega|)^{2s}|V_gf(x,\omega)|^2\,dx\,d\omega<\infty.
\end{equation}
We fix $x',\omega'\in\mathbb R^d$. By Proposition~\ref{est} and Lemma~\ref{adm1}, there exists $N \geq 0$ such that

\begin{eqnarray*}
&&\int_{\mathbb R^{3d}} |K_-(t,\omega')|\,|(\tau_{x'}g)(t)| |W_g(t)^{-1}||(\tau_xg)(t)|\,|K_-(\omega,t)|
|V_gf(x,\omega)|\,dt\,dx\,d\omega\\
&&\leq C_{d,g,\omega'} \int_{\mathbb R^d} (1+|t|)^{2\lambda+N} |(\tau_{x'}g)(t)| \left( \int_{\mathbb R^{2d}} (1+|\omega|)^\lambda |(\tau_xg)(t)||V_gf(x,\omega)|
\,dx\,d\omega \right)dt.
\end{eqnarray*}

We now estimate the inner integral. By the Cauchy-Schwarz inequality, and Lemma~\ref{translationest} we get
\begin{eqnarray*}
\int_{\mathbb R^{2d}} (1+|\omega|)^\lambda |(\tau_xg)(t)||V_gf(x,\omega)| \,dx\,d\omega  \! \! \! \! \! \! \! &&\leq
\left( \int_{\mathbb R^{2d}} (1+|\omega|)^{-2s} |(\tau_xg)(t)|^2\,dx\,d\omega
\right)^{1/2}\\
&& \times
\left(\int_{\mathbb R^{2d}}(1+|\omega|)^{2(\lambda+s)}|V_gf(x,\omega)|^2\,dx\,d\omega\right)^{1/2}
\end{eqnarray*}

Since $s>d/2$, by Lemma~\ref{est11} and Tonelli's theorem, there exists $C_g>0$ such that
$$
\int_{\mathbb R^{2d}}(1+|\omega|)^{-2s}|(\tau_xg)(t)|^2\,dx\,d\omega=
\left(\int_{\mathbb R^d}(1+|\omega|)^{-2s}\,d\omega\right)
\left(\int_{\mathbb R^d}|(\tau_xg)(t)|^2\,dx\right)\leq C_{d,s}C_g<\infty,
$$
We choose $s>d/2$. By Proposition~\ref{est}, Lemma~\ref{est11}, and the Cauchy-Schwarz
inequality, we obtain
\begin{eqnarray*}
&&\int_{\mathbb R^{3d}} |K_-(t,\omega')|\,|(\tau_{x'}g)(t)| |W_g(t)^{-1}||(\tau_xg)(t)|
|K_-(\omega,t)||V_gf(x,\omega)| \,dt\,dx\,d\omega\\
&&\leq C_{d,g,s,\omega'} \left(
\int_{\mathbb R^d} (1+|t|)^{2\lambda+N}|(\tau_{x'}g)(t)|\,dt \right) \left( \int_{\mathbb R^{2d}} (1+|\omega|)^{2(\lambda+s)}
|V_gf(x,\omega)|^2 \,dx\,d\omega \right)^{1/2} <\infty.
\end{eqnarray*}
Indeed, the first factor is finite since $\tau_{x'}g\in\mathcal S(\mathbb R^d;\mathbb R_{0,d})$, see Proposition \ref{tranS}, while the second factor is finite by \eqref{conV}.

So by Fubini's theorem we get
$$ V_gf(x',\omega')= \int_{\mathbb R^{2d}} \Bigg[ \frac{1}{(2\pi)^{d}}
\int_{\mathbb R^d}K_-(t,\omega') (\tau_{x'} g)(t)W_g(t)^{-1} \overline{(\tau_x g)(t)} K_-(\omega,t)\,dt\Bigg]V_gf(x,\omega)
d\omega\,dx.
$$
This proves \eqref{kenrel1}. Now, we prove that the integral \eqref{kenrel1} is absolutely convergent. By Proposition \ref{est} we have
\begin{equation}
	\label{e0}
	|K_-(t,\omega') K_-(\omega,t)|\leq C_{\omega,\omega'}(1+|t|)^{2\lambda}.
\end{equation}
For fixed $x,x'\in\mathbb R^d$, the functions $\tau_xg$ and $\tau_{x'}g$ belong to $\mathcal S(\mathbb R^d;\mathbb R_{0,d})$, see Proposition \ref{tranS}. Hence, for every $M>0$, we have
\begin{equation}
	\label{e1}
	|(\tau_xg)(t)|\leq C_{M,x}(1+|t|)^{-M}, \qquad |(\tau_{x'}g)(t)|\leq C_{M,x'}(1+|t|)^{-M}
\end{equation}
Since by hypothesis the function $g$ is an admissible function  by Lemma~\ref{adm1}, \eqref{e0} and \eqref{e1} we have
$$\left| K_-(t,\omega') (\tau_{x'}g)(t) W_g(t)^{-1} \overline{(\tau_xg)(t)} K_-(\omega,t)\right| \leq	C(1+|t|)^{2\lambda+N-2M}.$$
Now, we choose $M$ such that $2M>2\lambda+d+N$, thus the integral in \eqref{kenrel1} is absolutely convergent.
\end{proof}

%\begin{remark}
%If $g$ is radial, by \eqref{radial} we have
%$$$ W_g(t)=\int_{\mathbb R^d}g(t-x)^2\,dx=\|g\|_2^2.$$
 %Therefore, the kernel in \eqref{kenrel1} reduces to the one associated with the CSTFT with radial window, see \cite[Thm. 8.3]{DeMartinoCSTFT}.
 %Thus, the weight $W_g$ measures the effect of the generalized Clifford translation on the window and describes precisely the additional feature arising when the radiality assumption is removed.
%\end{remark}

\section{Uncertainty principles of the CSTFT}

Uncertainty principles are fundamental results in harmonic analysis, expressing the impossibility of arbitrarily localizing a function and its Fourier transform simultaneously; see \cite{HJ,Foll}. In the context of time-frequency analysis, they provide quantitative restrictions on the concentration of a signal in both time and frequency.

In this section, we establish some uncertainty principles for the CSTFT with non-radial window functions. In the classical setting, uncertainty principles for the STFT have been extensively studied; see, for example, \cite{BJ, C, GK,W}.  The uncertainty principles established in this section are obtained by transferring known inequalities for the Clifford-Fourier transform to the CSTFT. Although the underlying argument is common to all the inequalities, the extension to non-radial windows is not immediate. Indeed, unlike the radial case, the generalized translation does not preserve the classical translation structure.  The uniform estimates established in Section 3 provide upper bounds for the quadratic form $Q_g$, while the admissibility condition introduced in Section 4 supplies the lower bounds needed for the corresponding weighted uncertainty principles. The resulting inequalities therefore involve additional weights and window-dependent constants that are absent in the radial setting.

\subsection{Pitt's inequality}
Pitt's inequality shows that the $L^2$ norm of the Fourier transform weighted by the singular factor $|\omega|^{-2\beta}$ can be controlled by the $L^2$ norm of the original function with the weight $|t|^{2\beta}$, see \cite{Beck}. In \cite[Theorem~3.1]{LF}, an analogue of Pitt's inequality was established for the Clifford Fourier transform. We recall this result in the following theorem.

\begin{theorem}[Pitt's inequality for the CFT]
\label{pitt}
Let $d\ge2$ be even, let $f\in\mathcal S(\mathbb R^d;\mathbb{R}_{0,d})$, and let $0\le\beta<d/2$.  Then there exists a constant $C_\beta>0$, depending only on $d$ and $\beta$, such that
$$ \int_{\mathbb R^d}|\omega|^{-2\beta} |\mathcal F_{-}f(\omega)|^2\,d\omega \le C_\beta \int_{\mathbb R^d}|t|^{2\beta}
|f(t)|^2\,dt .
$$
\end{theorem}

We now establish a counterpart of Pitt's inequality for the CSTFT with a non-radial window function.
\begin{theorem}
\label{pittgeneral}
Let $d\geq 2$ be even and set $\lambda=\frac{d-2}{2}$. Let $f,g\in\mathcal{S}(\mathbb{R}^d;\mathbb{R}_{0,d})$, and let $0\leq\beta<d/2$. Then there exists a constant $C_{\beta,d}>0$, depending only on $d$ and $\beta$, such that
\begin{equation}
\int_{\mathbb{R}^{2d}} |\omega|^{-2\beta} |V_gf(x,\omega)|^2 \,dx\,d\omega
\leq C_{\beta,d}\, \left( 	\int_{\mathbb{R}^d} (1+|\xi|)^{2\lambda} |\mathcal{F}_-g(\xi)|^2\,d\xi\right) \left(\int_{\mathbb{R}^d}
(1+|t|)^{2\lambda+2 \beta} |f(t)|^2\,dt \right).
\label{pittgeneral1}
\end{equation}
\end{theorem}

\begin{proof}
We apply the Pitt's inequality for the Clifford Fourier transform, see Theorem \ref{pitt}, to the function $(\tau_xg)(t)f(t)$, with $x$ fixed, and we obtain
$$
\int_{\mathbb{R}^d}|\omega|^{-2\beta}|V_gf(x,\omega)|^2\,d\omega
\leq C_\beta \int_{\mathbb{R}^d} |t|^{2\beta} |(\tau_xg)(t)f(t)|^2\,dt.
$$
By the inequality \eqref{estmou} and integrating with respect to $x$, we have
\begin{equation}
\label{pittstep1}
\int_{\mathbb{R}^{2d}}|\omega|^{-2\beta}|V_gf(x,\omega)|^2 dx\,d\omega
\leq  C_{\beta,d} \int_{\mathbb{R}^d} |t|^{2\beta}|f(t)|^2 \left(
\int_{\mathbb{R}^d}|(\tau_xg)(t)|^2\,dx\right)dt.
\end{equation}
Finally, the result follows by Lemma \ref{translationest}.

\end{proof}

%\begin{remark}
%If we consider the function $g$ to be radial in the abovre result we can simply apply Theorem \ref{pitt} and get the following estimate
%$$\int_{\mathbb{R}^{2d}}| \omega|^{-2 \beta} |V_g f(x,\omega)|^2dx\, d\omega \le C_\beta \, \lVert g\rVert_{L^2(\mathbb{R}^d)}^2\int_{\mathbb{R}^d}\lvert t\rvert^{2\beta}\lvert f(t)\rvert^2\, dt,$$
%\end{remark}

The previous Pitt inequality yields a simple quantitative uncertainty principle for the frequency concentration of the CSTFT.

\begin{proposition}
\label{cosPitt}
Let $d\geq2$ be even, and let $0<\beta<d/2$. Let $f$, $g\in\mathcal S(\mathbb R^d;\mathbb R_{0,d})$, with $f\neq0$, and suppose that $g$ is an admissible window.	Assume that, for some $0\leq\varepsilon<1$ and $R>0$,
\begin{equation}
\label{iponew}
\int_{\mathbb R^d} \int_{|\omega|\leq R} |V_gf(x,\omega)|^2 d \omega dx
\geq (1-\varepsilon) \|V_gf\|_{L^2(\mathbb R^{2d}; \mathbb{R}_{0,d})}^2.
\end{equation}
Then there exist constants $c_g>0$ and $N\geq0$ such that
\begin{equation}
\label{thesis}
R \geq\left[\frac{(1-\varepsilon)c_g
	\int_{\mathbb R^d}(1+|t|)^{-N}|f(t)|^2\,dt}{C_{\beta,d}\mathcal A_g\int_{\mathbb R^d} (1+|t|)^{2\lambda+2 \beta}|f(t)|^2\,dt}\right]^{\frac{1}{2\beta}},
\end{equation}
where $ \mathcal{A}_g$ has been introduced in \eqref{new}.
\end{proposition}

\begin{proof}
By Theorem~\ref{par} and Tonelli's theorem, whose application is justified by Proposition~\ref{SchwartzCSTFT}, we obtain
$$
\|V_gf\|_{L^2(\mathbb R^{2d}; \mathbb{R}_{0,d})}^2= \int_{\mathbb R^d} \int_{\mathbb R^d} |(\tau_xg)(t)f(t)|^2
\,dt\,dx= \int_{\mathbb R^d} \left( \int_{\mathbb R^d} |(\tau_xg)(t)f(t)|^2\,dx \right)dt.
$$
Since $g$ is admissible, taking $u=f(t)$ in Definition \ref{adm} there exist constants $c_g>0$ and $N\geq0$ such that
$$ \|V_gf\|_{L^2(\mathbb R^{2d}; \mathbb{R}_{0,d})}^2 \geq c_g
 \int_{\mathbb R^d} (1+|t|)^{-N} |f(t)|^2\,dt.
$$
By \eqref{iponew} and the above inequality we get
\begin{equation}
\label{star7}
\int_{\mathbb R^d} \int_{|\omega|\leq R} |V_gf(x,\omega)|^2 \,d\omega\,dx \geq (1-\varepsilon)c_g \int_{\mathbb R^d} (1+|t|)^{-N} |f(t)|^2\,dt
\end{equation}
On the other hand, if $|\omega|\leq R$, we have $1\leq R^{2\beta}|\omega|^{-2\beta}$, hence by Theorem~\ref{pittgeneral} we have
$$
\int_{\mathbb R^d} \int_{|\omega|\leq R} |V_gf(x,\omega)|^2 \,d\omega\,dx \leq  R^{2\beta} \int_{\mathbb R^{2d}} |\omega|^{-2\beta} |V_gf(x,\omega)|^2 \,dx\,d\omega \leq C_{\beta,d}R^{2\beta}\mathcal A_g \int_{\mathbb R^d} (1+|t|)^{2\lambda+2\beta} |f(t)|^2\,dt.
$$
Finally by combining  \eqref{star7} and the above inequality we get
$$ (1-\varepsilon)c_g \int_{\mathbb R^d} (1+|t|)^{-N} |f(t)|^2\,dt \leq
C_{\beta,d}R^{2\beta}\mathcal A_g \int_{\mathbb R^d} (1+|t|)^{2\lambda+2\beta} |f(t)|^2\,dt.$$
This implies \eqref{thesis}.
\end{proof}

\begin{remark}
The previous result gives a quantitative interpretation of an uncertainty principle in terms of the frequency concentration of the CSTFT. Indeed, the assumption \eqref{iponew} means that at least a fraction $1-\varepsilon$ of the total energy of the CSTFT is contained in the frequency ball $B(0,R)$. Proposition~\ref{cosPitt} shows that the radius $R$ of the ball cannot be arbitrarily small.
\end{remark}

\subsection{Logarithmic uncertainty principle for the CSTFT}

In classical harmonic analysis, Pitt's inequality is used to establish the so-called logarithmic uncertainty principle for the Fourier transform; see \cite{Beck}. This principle provides a quantitative description of the simultaneous concentration of a function and its Fourier transform in terms of logarithmic moments in the space and frequency variables. An analogue of this result for the Clifford-Fourier transform was established in \cite{LF}.

\begin{theorem}
\label{thmlog}
Let $d\geq2$ be even and $f\in\mathcal S(\mathbb R^d;\mathbb{R}_{0,d})$.  Then
$$ \int_{\mathbb R^d}\log|\omega|\, |\mathcal F_{-}f(\omega)|^2\,d\omega
+\int_{\mathbb R^d}\log |t|\,|f(t)|^2\,dt \ge \gamma_d\|f \|_{L^2 (\mathbb{R}^d; \mathbb{R}_{0,d})}^2,
$$
where $\gamma_d:=\left(\psi\!\left(\frac d4\right)+\ln 2\right)$ and  $\psi(t)= \frac{\Gamma'(t)}{\Gamma(t)}$.
\end{theorem}

We now turn our attention to a logarithmic uncertainty principle for the CSTFT.
\begin{theorem}[Logarithmic uncertainty principle for non-radial windows]
\label{lognonradial}
Let $d\geq2$ be even and let $f\in\mathcal S(\mathbb R^d;\mathbb R_{0,d})$. If $d=2$, let
$g\in\mathcal S(\mathbb R^2;\mathbb R_{0,2})$. Then
\begin{equation}
\label{loddime2}
\! \! \! \! \! \!
\int_{\mathbb R^4} \log|\omega||V_gf(x,\omega)|^2 dx d\omega+
\|g\|_{L^2(\mathbb R^2;\mathbb R_{0,2})}^2\int_{\mathbb R^2} \log|t| |f(t)|^2 dt \geq
\gamma_2\|g\|_{L^2(\mathbb R^2;\mathbb R_{0,2})}^2 \|f\|_{L^2(\mathbb R^2;\mathbb R_{0,2})}^2.
\end{equation}
If $d\geq4$, assume that $g\in\mathcal S(\mathbb R^d;\mathbb R_{0,d})$ is an admissible
window. Then there exist $c_g>0$ and $N \geq 0$, such that
\begin{align}
&\int_{\mathbb R^{2d}}\log|\omega|\,|V_gf(x,\omega)|^2\,dx\,d\omega+ C_d\mathcal{A}_g\int_{|t|\geq1}  \log|t|\,(1+|t|)^{d-2} |f(t)|^2\,dt
\nonumber\\
&+c_g\int_{|t|<1} \frac{\log|t|}{(1+|t|)^N}  |f(t)|^2\,dt
\geq \gamma_d c_g \int_{\mathbb R^d} \frac{|f(t)|^2} {(1+|t|)^N}\,dt .
\label{logshort}
\end{align}
where the constant $\gamma_d$ is defined in Theorem \ref{thmlog}, $ \mathcal{A}_g$ is defined in \eqref{new} and $C_d$ is a constant depending on the dimension.
\end{theorem}

\begin{proof}
For fixed $x\in\mathbb R^d$, by hypothesis and Proposition \ref{tranS} we have that $\tau_xg \cdot f\in\mathcal S(\mathbb R^d;\mathbb R_{0,d})$.  By applying Theorem \ref{thmlog} to $\tau_xg \cdot f$ we have
$$\int_{\mathbb R^d} \log|\omega|\, |V_gf(x,\omega)|^2\,d\omega +
\int_{\mathbb R^d} \log|t|\, |(\tau_xg)(t)f(t)|^2\,dt \geq \gamma_d \int_{\mathbb R^d} |(\tau_xg)(t)f(t)|^2\,dt.$$
Integrating with respect to $x$ and using Fubini's theorem we get
\begin{equation}
\label{logW}
\int_{\mathbb R^{2d}} \log|\omega|\, |V_gf(x,\omega)|^2 \,dx\,d\omega+
\int_{\mathbb R^d} \log|t|\,Q_g(t, f(t))\,dt \geq \gamma_d
\int_{\mathbb R^d} Q_g(t, f(t))\,dt,
\end{equation}
where $Q_g(t, f(t))$ has been introduced in \eqref{Qinte}. The inequality \eqref{loddime2} follows by \eqref{d2}. The application of Fubini's theorem is justified by the fact that both $V_gf(x,\omega)$ and $(\tau_xg)(t)f(t)$ belong to $\mathcal{S}(\mathbb R^{2d};\mathbb R_{0,d})$ (see Propositions~\ref{SchwartzCSTFT} and~\ref{tranS}), while the logarithmic weights are locally integrable and have at most polynomial growth at infinity.
\\Now, we focus on the case $d \geq 4$. 
If $|t|\geq1$, then $\log|t|\geq0$. By \eqref{Qgtwosided} we get
$$ \log|t|\,Q_g(t,f(t)) \leq C_d \log|t|  \mathcal{A}_g (1+|t|)^{d-2} |f(t)|^2.$$
If $|t|<1$, then $\log|t|<0$. Since $g$ is admissible there exist $c_g>0$ and $N \geq 0$ such that we have
$$ \log|t|\,Q_g(t,f(t)) \leq c_g \frac{\log|t|}{(1+|t|)^N} |f(t)|^2.$$
We estimate the logarithmic term in \eqref{logW} as
\begin{equation}
\label{logWl}
\! \! \! \! \! \! \! \!
\int_{\mathbb R^d} \log|t|\,Q_g(t,f(t))\,dt \leq  C_d \mathcal{A}_g\int_{|t|\geq1} \log|t|\,(1+|t|)^{d-2}
|f(t)|^2\,dt+c_g \int_{|t|<1} \frac{\log|t|}{(1+|t|)^N} |f(t)|^2\,dt.
\end{equation}
Finally, for even $d\geq4$ the constant given by $\gamma_d =\psi\left(\frac d4\right)+\log2$ is positive. Since $g$ is admissible there exist $c_g>0$ and $N \geq 0$ such that we have
\begin{equation}
\label{Qlower} \gamma_d\int_{\mathbb R^d}Q_g(t,f(t))\,dt\geq \gamma_d c_g \int_{\mathbb R^d}
\frac{|f(t)|^2}{(1+|t|)^N}\,dt.
\end{equation}
Combining \eqref{logW}, \eqref{logWl}, and \eqref{Qlower} we get \eqref{logshort}.
\end{proof}

\begin{remark}
If $g\in\mathcal S(\mathbb R^d)$ is a nonzero real-valued radial function, then $Q_g(t,u)=\|g\|_2^2|u|^2$. Consequently, inequality \eqref{logW} yields the following logarithmic uncertainty principle for the CSTFT:
$$\int_{\mathbb R^{2d}}\log|\omega||V_gf(x,\omega)|^2\,dx\,d\omega+\|g\|_{L^2(\mathbb R^d)}^2
\int_{\mathbb R^d}\log|t|\,|f(t)|^2\,dt\geq \gamma_d\|g\|_{L^2(\mathbb R^d)}^2 \|f\|_{L^2(\mathbb R^d;\mathbb R_{0,d})}^2.$$
\end{remark}

\begin{remark}
The distinction between the cases $d=2$ and $d>2$ reflects a structural feature of the CSTFT transform. In dimension two, the generalized translation coincides with the classical translation and so no additional condition on the window is required.
For even dimensions $d>2$, the generalized translation does not coincide, in general, with the classical one. This motivates the introduction of admissible windows to ensure the invertibility of the associated weight.
\end{remark}

\subsection{Heisenberg-Pauli-Weyl uncertainty principle for the CSTFT}

In the classical setting, the logarithmic uncertainty principle implies the well-known Heisenberg-Pauli-Weyl uncertainty principle; see \cite{Hein}. We now present the counterpart of the Heisenberg uncertainty principle for the Clifford-Fourier transform.

\begin{theorem}
\label{HeinCF}
Let $d\ge2$ be even and let $f\in\mathcal S(\mathbb R^d;\mathbb{R}_{0,d})$ be nonzero.  Then
$$ \left( \int_{\mathbb R^d}|t|^2|f(t)|^2\,dt\right)^{1/2} \left( \int_{\mathbb R^d}|\omega|^2
|\mathcal F_{-}f(\omega)|^2\,d\omega \right)^{1/2} \ge e^{\gamma_d}\|f\|_2^2,
$$
with $\gamma_d$ as in Theorem~\ref{thmlog}. 
\end{theorem}

We now turn our attention to the Heisenberg uncertainty principle for the CSTFT. In this setting, the aim is to establish a lower bound involving the second moment of the signal in the time variable and the second moment of its CSTFT in the frequency variable.

\begin{theorem}[Heisenberg uncertainty principle for the CSTFT]
\label{thm:Heisenberg-CSTFT}
Let $d\geq2$ be even and let $f,g\in\mathcal S(\mathbb R^d;\mathbb R_{0,d})$, with $f\neq0$ and $g\neq0$. If $d=2$, then
\begin{equation}
\label{Hd2}
\left(\int_{\mathbb R^2}|t|^2|f(t)|^2\,dt\right)^{1/2} \left( \int_{\mathbb R^4} |\omega|^2|V_gf(x,\omega)|^2\,dx\,d\omega \right)^{1/2} \geq e^{\gamma_2}\|g\|_{L^2(\mathbb R^2;\mathbb R_{0,2})} \|f\|_{L^2(\mathbb R^2;\mathbb R_{0,2})}^2.
\end{equation}
If $d\geq4$ and $g$ is an admissible window, then there exist $c_g>0$ and $N \geq 0$ such that
\begin{equation}
\label{Hstft}
\left( \int_{\mathbb R^d} (1+|t|)^{d}|f(t)|^2\,dt \right)^{1/2}\left(
\int_{\mathbb R^{2d}} |\omega|^2|V_gf(x,\omega)|^2\,dx\,d\omega \right)^{1/2} \geq
\frac{c_d c_g}{\sqrt{\mathcal{A}_g}} \int_{\mathbb R^d} \frac{|f(t)|^2}{(1+|t|)^N}\,dt,
\end{equation}
where $ \mathcal{A}_g$ has been introduced in \eqref{new} and $c_d$ is a constant depending on the dimension.
\end{theorem}

\begin{proof}
For fixed $x\in\mathbb R^d$, and by Proposition \ref{tranS} we have $\tau_xg \cdot f \in\mathcal S(\mathbb R^d;\mathbb R_{0,d})$, thus the Heisenberg uncertainty principle for the Clifford-Fourier transform, see Theorem \ref{HeinCF}, gives

\begin{equation}
\label{H1}
\left( \int_{\mathbb R^d} |t|^2|(\tau_xg)(t)f(t)|^2\,dt \right)^{1/2}
\left( \int_{\mathbb R^d} |\omega|^2|V_gf(x,\omega)|^2\,d\omega\right)^{1/2}
\geq c_d \int_{\mathbb R^d} |(\tau_xg)(t)f(t)|^2\,dt.
\end{equation}
Integrating \eqref{H1} with respect to $x$ and applying the Cauchy-Schwarz inequality and Tonelli's theorem we get

\begin{equation}
\label{H2}
\left(\int_{\mathbb R^d} |t|^2Q_g(t,f(t))\,dt \right)^{1/2}
\left( \int_{\mathbb R^{2d}} |\omega|^2|V_gf(x,\omega)|^2\,dx\,d\omega \right)^{1/2} \geq
c_d \int_{\mathbb R^d}Q_g(t,f(t))\,dt,
\end{equation}
where $Q_g(t,f(t))$ has been defined in \eqref{Qinte}. Thus the inequality \eqref{Hd2} follows by \eqref{d2}. Let now $d\geq4$. By the upper estimate of \eqref{Qgtwosided}, we get
$$ \int_{\mathbb R^d}|t|^2Q_g(t,f(t))\,dt \leq c_d \mathcal{A}_g \int_{\mathbb R^d} |t|^2(1+|t|)^{d-2}|f(t)|^2\,dt \leq c_d \mathcal{A}_g \int_{\mathbb R^d} (1+|t|)^{d}|f(t)|^2\,dt,$$
whereas by using the hypothesis that $g$ is admissible there exists $c_g>0$ and $N \geq 0$ such that
$$ \int_{\mathbb R^d}Q_g(t,f(t))\,dt \geq c_g \int_{\mathbb R^d} \frac{|f(t)|^2}{(1+|t|)^N}\,dt.
$$	
Using these two estimates in \eqref{H2}, we obtain \eqref{Hstft}.
\end{proof}

\begin{remark}
If $g\in\mathcal S(\mathbb R^d)$ is a nonzero real-valued radial function, then $Q_g(t,u)=\|g\|_2^2|u|^2$. Consequently, inequality~\eqref{H2} yields
$$ \left(\int_{\mathbb R^{2d}}|\omega|^2 |V_gf(x,\omega)|^2 dx\,d\omega\right)^{1/2}\left(\int_{\mathbb R^d}|t|^2|f(t)|^2\,dt\right)^{1/2} \geq c_d\,\|g\|_{L^2(\mathbb R^d)} \|f\|_{L^2(\mathbb R^d;\mathbb R_{0,d})}^2.
$$
\end{remark}

\begin{remark}
The uncertainty principles obtained above share a common underlying mechanism, but their non-radial counterparts exhibit a different structure from the classical ones. In particular,
the quadratic form $Q_g(t,u)$ naturally replaces the constant window norm appearing in the radial case. The resulting weighted inequalities reflect the nontrivial behavior of generalized translations and recover the classical form when $g$ is a real-valued radial window.
\end{remark}

\section{Concluding remarks}

In this paper, we have extended several fundamental properties of the Clifford short-time Fourier transform to non-radial window functions in even dimensions. The main difficulty arises from the generalized translation operator, which no longer coincides with
the classical translation. By establishing an angular estimate for the Clifford-Fourier kernel, we obtained a uniform $L^2$-bound for generalized translates of
Schwartz functions. This allowed us to derive a weighted Moyal identity and, under an admissibility condition on the window function, reconstruction and reproducing kernel formulas. We also established Pitt-type, logarithmic, and Heisenberg-type uncertainty
principles in this setting.

Several questions deserve further investigation. In particular, it would be interesting to establish Hardy, Beurling, and Cowling-Price-type uncertainty principles for the CSTFT with non-radial windows. Another natural direction is the extension of the present results to odd dimensions, where the absence of an explicit closed-form expression for the Clifford-Fourier kernel presents additional difficulties; see~\cite{Lian}.

\section*{Acknowledgements and AI Declaration}
The second author acknowledges the support of the China Scholarship Council (CSC) under Grant No. 202506030086.
During the preparation of this work the authors used ChatGPT in order to check some computations. After using this tool/service, the authors reviewed and edited the content as needed and take full responsibility for the content of the published article.


\begin{thebibliography}{99}
		
		\bibitem{AS} M.~Abramowitz, I.A.~Stegun, \emph{Handbook of mathematical functions with formulas, graphs, and mathematical tables}, National Bureau of Standards Applied Mathematics Series, No. 55, U. S. Government Printing Office, Washington, DC. 
		
\bibitem{AB} L.D.~Abreu, \emph{Sampling and interpolation in Bargmann-Fock spaces of polyanalytic functions}, Appl. Comput. Harmon. Anal., 29, (2010), 287-302 .


		
\bibitem{ADDS} D.~Alpay, A.~De Martino, K.~Diki, D.C.~Struppa, \emph{Short-time Fourier transform and superoscillations}, Appl. Comput. Harmon. Anal.  73 (2024), Paper No. 101689, 31 pp.
		
\bibitem{Beck} W.~Beckner, \emph{Pitt’s inequality and the uncertainty principle}, Proc. Am. Math. Soc. 123, 1897–1905, 1995.
		
	\bibitem{BJ} A.~Bonami, B.~Demange, P.~Jaming, \emph{Hermite functions and uncertainty principles for the Fourier and the windowed Fourier transforms}, Rev. Mat. Iberoamericana {\bf 19} (2003), no.~1, 23-55.
		
	\bibitem{red} F.~Brackx, R.~Delanghe, F.~Sommen, {\em Clifford Analysis}, Pitman Res. Notes in Math., 76, 1982.
		
	\bibitem{BSS} F.~Brackx, N.~De Schepper, F.~Sommen, \emph{The Clifford-Fourier transform}, J. Fourier Anal. Appl. {\bf 11} (2005), no.~6, 669-681.
	
		\bibitem{BSS1}
	F.~Brackx, N.~De Schepper, F.~Sommen, \emph{The two-dimensional Clifford-Fourier transform}, J. Math. Imaging Vision {\bf 26} (2006), no.~1-2, 5-18.
	

	
	\bibitem{C} L.~Cohen, \emph{The uncertainty principle for the short-time Fourier transform and wavelet transform}, in {\it Wavelet transforms and time-frequency signal analysis}, 217-232, Appl. Numer. Harmon. Anal., Birkh\"auser Boston, Boston, MA.
	
		\bibitem{ColomboSabadiniStruppa2011}
	F. Colombo, I. Sabadini, and D. C. Struppa, \emph{Noncommutative Functional Calculus. Theory and Applications of Slice Hyperholomorphic Functions}, Birkh\"auser, Basel, 2011.


\bibitem{DeBieCFT} H.~De~Bie, Y.~Xu, \emph{On the Clifford-Fourier transform}, Int. Math. Res. Not. IMRN 2011 (2011), no.~22, 5123-5163.

\bibitem{DeBieFourier} H.~De~Bie, \emph{Fourier transforms in Clifford analysis},
in Operator Theory, Springer Basel, Basel, 2015, 1651-1672.

\bibitem{DeJ} M.F.E.~De Jeu,\emph{The Dunkl transform}. Invent. Math. 113 (1993), 147-162.

\bibitem{DK}A.~De Martino, K.~Diki, \emph{On the quaternionic short-time Fourier and Segal-Bargmann transforms}, Mediterr. J. Math. {\bf 18} (2021), no.~3, Paper No. 110, 22 pp.

\bibitem{DK1} A.~De Martino, K.~Diki, \emph{On the polyanalytic short-time Fourier transform in the quaternionic setting}, Commun. Pure Appl. Anal. {\bf 21} (2022), no.~11, 3629-3665.


\bibitem{DeMartinoCSTFT} A.~De~Martino, \emph{On the Clifford short-time Fourier transform and its properties}, Appl. Math. Comput. 418 (2022), Article No.~126812.

				
\bibitem{green} R.~Delanghe, F.~Sommen, V.~Soucek, \emph{Clifford algebra and spinor-valued functions. A function theory for the Dirac operator}, volume 53 of mathematics and its Applications. Kluwer Academic Publishers Group, Dordrecht (1992).


\bibitem{Don} D.L.~Donoho, P.B.~Stark, \emph{Uncertainty principles and signal recovery}, SIAM J. Appl. Math. 49 (1989), no. 3, 906–931.

\bibitem{Dunkl} C.F.~Dunkl, \emph{Hankel transforms associated to finite reflection groups, in Proc. of Special Session on Hypergeometric Functions on Domains of Positivity}, Jack Polynomials and Applications (Tampa, 1991), Contemp.Math. 138 (1992), 123-138.

\bibitem{GhobberJaming}
S.~Ghobber, P.~Jaming, \emph{Uncertainty principles for integral operators}, Studia Math. {\bf 220} (2014), no.~3, 197-220. 

\bibitem{G1} K.~Gröchenig, \emph{Uncertainty Principles for Time-Frequency Representations}. In: Feichtinger, H.G., Strohmer, T. (eds) Advances in Gabor Analysis. Applied and Numerical Harmonic Analysis. Birkhäuser, Boston, MA. (2003).

\bibitem{G}K.~Gr\"ochenig, \emph{Foundations of time-frequency analysis}, Applied and Numerical Harmonic Analysis, Birkh\"auser Boston, Boston, MA, 2001.


\bibitem{GK}K.~Gröchenig, G. Zimmermann, \emph{Hardy’s theorem and the short-time Fourier transform of Schwartz functions}, J. London Math. Soc. 63 (2001), 205–214.



\bibitem{Hein} W.~Heisenberg, \emph{Über den anschaulichen Inhalt der quantentheoretischen Kinematik und Mechanik}, Z. Physik 43, 172–198, 1927.



\bibitem{H} R.~Howe, \emph{The mathematical heritage of Hermann Weyl}, Proceedings of Symposia in Pure Mathematics, 48, American Mathematical Society, Providence, RI, 1988.

\bibitem{HJ} V.P.~Khavin, B. J\"oricke, {\it The uncertainty principle in harmonic analysis}, Ergebnisse der Mathematik und ihrer Grenzgebiete (3), 28, Springer, Berlin, 1994.


\bibitem{Foll} G.B.~Folland, A.~Sitaram, \emph{The uncertainty principle: a mathematical survey}, J. Fourier Anal. Appl. 3 (1997), no. 3, 207–238.

\bibitem{LiLengFeiPW} S.~Li, J.~Leng, M.~Fei, \emph{Paley-Wiener-type theorems for the Clifford-Fourier transform}, Math. Methods Appl. Sci. 42 (2019), no.~18, 6101-6113.

\bibitem{LF} S.~Li, M.~Fei, \emph{Pitt's inequality and logarithmic uncertainty principle for the Clifford-Fourier transform}, Adv. Appl. Clifford Algebras 33 (2023), no.~1, Paper No.~2.

\bibitem{Lian} P.~Lian, \emph{The bounds of the odd dimensional Clifford-Fourier kernels}, Ann. Mat. Pura Appl. (4) {\bf 201} (2022), no.~3, 1213-1228.

\bibitem{lieb} E.H.~Lieb, \emph{Integral bounds for radar ambiguity functions and Wigner distributions}, J. Math. Phys. {\bf 31} (1990), no.~3, 594-599.





\bibitem{T} S.~Thangavelu, \emph{A scalar-valued Fourier transform for the Heisenberg group, From classical analysis to analysis on fractals}, vol. 1, Appl. Numer. Harmon. Anal., Birkhäuser/Springer, Cham, 2023, pp. 137–163. 

\bibitem{W} E.~Wilczok, \emph{New uncertainty principles for the continuous Gabor transform and the continuous wavelet transform}, Doc. Math. {\bf 5} (2000), 201-226.





	\end{thebibliography}
\end{document}